\documentstyle{amsppt}
\magnification 1000 \hsize=6.1in \vsize=8.75in
\input amstex
\TagsOnRight

\document
\topmatter

\title
Convolution Dirichlet series and a Kronecker limit formula for
second-order Eisenstein series
\endtitle  
\author  Jay Jorgenson and Cormac O'Sullivan \endauthor
\date 21 March 2004\enddate
\endtopmatter
\NoRunningHeads


\def\H{{\bold H}}
\def\F{{\frak F}}
\def\C{{\Bbb C}}
\def\R{{\Bbb R}}
\def\Z{{\Bbb Z}}

\def\N{{\Bbb N}}
\def\G{{\Gamma}}
\def\GH{{\G \backslash \H}}
\def\g{{\gamma}}

\def\ee{{\varepsilon}}
\def\K{{\Cal K}}
\def\Re{\text{\rm Re}}
\def\Im{\text{\rm Im}}
\def\PSL{\text{\rm PSL}}

\def\ca{{\frak a}}
\def\cb{{\frak b}}
\def\cc{{\frak c}}

\def\ci{{\infty}}

\def\sa{{\sigma_\frak a}}
\def\sb{{\sigma_\frak b}}
\def\sc{{\sigma_\frak c}}

\def\p{\endproclaim \flushpar {\bf Proof: }}

\def\bs{\ \ $\qed$ \vskip 3mm}


\head \bf Abstract\rm
\endhead

In this article we derive analytic and Fourier aspects of a
Kronecker limit formula for second-order Eisenstein series. Let
$\Gamma$ be any Fuchsian group of the first kind which acts on the
hyperbolic upper half-space $\H$ such that the quotient $\Gamma
\backslash \H$ has finite volume yet is non-compact. Associated to
each cusp of $\Gamma \backslash \H$, there is a classically
studied \it first-order \rm non-holomorphic Eisenstein series
$E(s,z)$ which is defined by a generalized Dirichlet series that
converges for $\Re (s) > 1$. The Eisenstein series $E(s,z)$ admits
a meromorphic continuation with a simple pole at $s=1$.
Classically, Kronecker's limit formula is the study of the
constant term $\K_{1}(z)$ in the Laurent expansion of $E(s,z)$ at
$s=1$. A number of authors recently have studied what is known as
the \it second-order \rm Eisenstein series $E^{\ast}(s,z)$, which
is formed by twisting the Dirichlet series that defines the
 series $E(s,z)$ by periods of a given
cusp form $f$. In the work we present here, we study an analogue
of Kronecker's limit formula in the setting of the second-order
Eisenstein series $E^{\ast}(s,z)$, meaning we determine the
constant term $\K_{2}(z)$ in the Laurent expansion of
$E^{\ast}(s,z)$ at its first pole, which is also at $s=1$.  To
begin our investigation, we prove a bound for the Fourier
coefficients associated to the first-order Kronecker limit
function $\K_{1}$. We then define two families of convolution
Dirichlet series, denoted by $L^{+}_{m}$ and $L^{-}_{m}$ with $m
\in \N$, which are formed by using the Fourier coefficients of
$\K_{1}$ and the weight two cusp form $f$.  We prove that for all
$m$, $L^{+}_{m}$ and $L^{-}_{m}$ admit a meromorphic continuation
and are holomorphic at $s=1$.  Turning our attention to the
second-order Kronecker limit function $\K_{2}$,  we first express
$\K_{2}$ as a solution to various differential equations.  Then we
obtain its complete Fourier expansion in terms of the cusp form
$f$, the Fourier coefficients of the first-order Kronecker limit
function $\K_{1}$, and special values $L^{+}_{m}(1)$ and
$L^{-}_{m}(1)$ of the convolution Dirichlet series. Finally, we
prove a bound for the special values $L^{+}_{m}(1)$ and
$L^{-}_{m}(1)$ which then implies a bound for the Fourier
coefficients of $\K_{2}$.  Our analysis leads to certain natural
questions concerning the holomorphic projection operator, and we
conclude this paper by examining certain numerical examples and
posing questions for future study.

\head \bf Table of Contents\rm
\endhead
\roster
\item"{1.}" Introduction and statement of results
\item"{2.}" The Fourier expansion of $\K_{2}$
\item"{3.}" Poincar\'e series and holomorphic projection
\item"{4.}" $\K_{2}$ and the holomorphic projection of $f\K_{1}$
\item"{5.}" The Dirichlet series $L_{m}^{+}$ and the holomorphic projection of $f\K_{1}$
\item"{6.}" The Dirichlet series $L_{m}^{-}$ and the proofs of Theorems 1.4 and 1.5
\item"{7.}" Bounding the Fourier coefficients of $\K_{1}$ and $\K_{2}$
\item"{8.}" Poincar\'e series: Proofs of Theorems 3.1 and 3.2
\item"{9.}" Proofs of Proposition 3.3 and the meromorphic continuation of $L_{m}^{+}$ and $L_{m}^{-}$
\item"{10.}" Examples
\item"{}" References
\endroster

\vskip .25in \head \bf \S 1. Introduction and statement of results
\rm
\endhead

\vskip .10in Let $\G$ contained in $\PSL_2(\R)$ be a Fuchsian
group of the first kind acting on the upper half plane $\H$ with
non-compact quotient $\G\backslash\H$. As usual, we write
$x+iy=z\in \H$. Set $V$ equal to the hyperbolic volume of
$\G\backslash\H$.  Assuming there is a cusp at $\infty$, let
$\G_\ci=\{ \g \in \G \,| \, \g \ci = \ci\}$, and, for simplicity
we may assume that $\G_\ci$ is generated by $z\mapsto z+1$. The
\it first-order \rm non-holomorphic Eisenstein series is defined
by the series
$$
E(z,s)=\sum_{\g \in \G_\infty\backslash\G} \Im(\g z)^s
$$
which converges for $\Re (s)>1$ and has a meromorphic continuation
to all $s$ in $\C$ (see, for example, Chapter 6 of  \cite{Iw1}).
The function $E(z,s)$ is known to have a simple pole at $s=1$ with
residue $V^{-1}$, so then, when denoting the constant part at
$s=1$ by $\K_1$, we can write
$$
E(z,s)=\frac{V^{-1}}{s-1}+\K_1(z)+O(s-1) \,\,\,\,\,\text{\rm as $s
\rightarrow 1$.}
$$
The first result which is known as \it Kronecker's first limit
formula \rm is the following.  If $\G=\PSL_2(\Z)$, then
$$
\K_1(z)=\frac{-1}{4\pi}\log(y^{12}|\Delta(z)|^2)+\frac{3}{\pi}(\g-\log
4\pi)\tag 1.1
$$
where
$$
\Delta(z) = e^{2\pi i z}
\prod\limits_{n=1}^{\infty}\left(1-e^{2\pi i n z}\right)^{24}
$$
is the discriminant function,
a weight 12 holomorphic cusp form for $\PSL_2(\Z)$, and $\g$ is
Euler's constant. Kronecker's second limit formula, which for
brevity we do not state here, is a determination of the constant
term at the first pole of the first-order non-holomorphic
Eisenstein series obtained by twisting the series definition of
$E(z,s)$ with a unitary character of $\G$.  We refer to \cite{La},
\cite{Si}, or \cite{Za1} for proofs of these classical results.

\vskip .10in Many generalizations of the Kronecker limit formulas
exist, and the results have diverse applications.  In \cite{La},
\cite{Si}, and \cite{Za2}, formulas for class numbers of algebraic
number fields are obtained; in \cite{C-P} and \cite {P-W}, the
limit formulas are used to find values of $|\eta(z)|$ at
quadratic irrationalities; in \cite{B-C-Z} and \cite{Ra}, special
values of the Rogers-Ramanujan continued fraction are evaluated;
and in \cite{R-S}, the limit formulas are used to explicitly
evaluate analytic torsion for flat line bundles on elliptic curves.
The analogue of Kronecker's first limit formula to Hilbert modular
varieties has been studied, beginning with \cite{As} for totally
real fields and \cite{E-G-M} for imaginary quadratic fields, then
\cite{J-L} for general number fields. Returning to the setting of
$\PSL_2(\R)$, the limit function $\K_1(z)$ has been determined for
other groups in \cite{Gn}; specific results for the Hecke
congruence subgroups $\G_0(N)$ are given in section 10 below.

\vskip .10in Our focus in this paper is to find formulas for the
constant part at $s=1$ of second-order Eisenstein series, which
are defined by twisting the classical non-holomorphic Eisenstein
series by a modular symbol.  In general, a non-holomorphic
second-order Eisenstein series $E^{\ast}(z,s)$ is associated to the following 
data: A Fuchsian group $\G$ of the first kind; a
parabolic subgroup of $\G$; and a weight two holomorphic form
which vanishes in each cusp of $\G$. The precise definition is given 
below. The series $E^{\ast}(z,s)$
was first defined and studied in \cite{Gd} in order to provide
another approach to the ABC-conjecture, which itself is connected
to a number of fundamental and motivating problems in number
theory, such as: Mordell's conjecture (a theorem of Faltings);
Szpiro's conjecture; the degree conjecture; Goldfeld's period
conjecture; and various questions and assertions regarding the
Shafarevich-Tate group.  In particular, we refer the reader to
\cite{Gd2} where Goldfeld states what he calls the Modular Symbol
Conjecture, together with a summary of the inter-relations between
the aforementioned conjectures as well as the role played by the
Modular Symbol Conjecture.  In \cite{M-M}, Manin and Marcolli
generalized the classical Gauss-Kuzmin theorem having to do with
the distribution of continued fractions.  Going further, the
authors 
develop connections between weighted averages of modular symbols, such as $E^*(z,s)$, and the distribution of continued fractions. The distribution of modular symbols themselves is elaborated by Petridis and Risager in \cite{P-R} with their work on $E^*(z,s)$ and its generalizations.
  In
\cite{K-Z}, Kleban and Zagier studied crossing probabilities and
free energies for conformally invariant critical 2-D systems,
which they derive from conformal field theory and certain
stochastic integrals.  It is shown in \cite{K-Z} that the crossing
probabilities and partition functions they encountered may be
expressed as values of what should now 
be viewed as holomorphic second-order modular forms.  As discussed
in the concluding remarks of \cite{K-Z}, second-order forms in
general can, in certain cases, be viewed as components of vector-valued modular forms
associated to certain representations of the Fuchsian group $\G$
into $\text{\rm SL}_{2}$.  In this way, the non-holomorphic
second-order Eisenstein series, and second-order forms in general,
are manifest in all aspects of the spectral theory, holomorphic
function theory, number theory, and algebraic geometry of certain
vector-valued functions on Riemann surfaces.  In summary,
second-order forms, which include $E^{\ast}(z,s)$, have at the
present an established place in number theory  \cite{Gd2}, \cite{M-M}, \cite{P-R} and in physics \cite{K-Z};
furthermore, additional connections to converse theorems in number theory, to spectral theory and
to algebraic geometry are pending.  As a result, any and all results
regarding second-order forms should be viewed as interesting for
their own sake as well as having wide yet unforeseen consequences.
 
\vskip .10in For the purposes of narrowing our attention, we will
concentrate on two aspects of the Kronecker limit formula:
Differential equations, and Fourier expansions, with the latter
necessarily requiring the study of the growth of the Fourier
coefficients. Before stating our results, let us establish
necessary background material and notation.

\vskip .10in Let $S_k(\G)$ be the space of holomorphic weight $k$
cusp forms for $\G$, meaning the vector space of holomorphic
functions $g$ on $\H$ which satisfy the transformation property
$$
g(\g z)=j(\g,z)^k g(z) \text{ \ \ with \ \ }j(\left(\smallmatrix a
& b \\ c & d \endsmallmatrix\right),z)=cz+d \text{ \ \ for \ \
}(\smallmatrix a & b \\ c & d \endsmallmatrix)\in \G,
$$
and decay rapidly in each cusp in the quotient space
$\G\backslash\H$.  As usual, we equip the vector space $S_{k}(\G)$
with the well-known Petersson inner product. Since the analytic
transformation $z \mapsto z + 1$ corresponds to an element of
$\G$, we have that any $f \in S_{k}(\G)$ admits a Fourier
expansion, for which we use the notation
$$
f(z)=\sum_{n=1}^\infty a_ne(nz) \,\,\,\,\, \text{\rm where
\,\,\,\,\,$e(z) = e^{2\pi i z}$,}
$$
and from which we define
$$
F(z)=\sum_{n=1}^\infty \frac{a_n}{n}e(nz)=2\pi
i\int_{i\infty}^{z}f(w)\,dw.
$$
For the remainder of this paper we set $f$ to have weight two: $f \in S_{2}(\G)$.
The \it modular symbol \rm $\langle \cdot, f \rangle $ associated
to $f$ is the homomorphism from $\G$ to $\C$ given by
$$
\langle\g,f\rangle=2\pi i\int_z^{\g z}f(w)\,dw=F(\g z)-F(z).
$$
The \it second-order \rm non-holomorphic Eisenstein series
associated to $f$ is defined for $\Re (s) >1$ by the convergent
$$
E^*(z,s)=\sum_{\g \in \G_\infty\backslash\G}
\langle\g,f\rangle\Im(\g z)^s.
$$
For any $\gamma, \tau \in \G$, the non-holomorphic Eisenstein
series satisfy the transformation properties
$$
\align
&E(\g z,s)-E(z,s)=0,\tag1.2\\
E^*(\g\tau z,s)-&E^*(\g z,s) -E^*(\tau z,s)+E^*(z,s)=0.\tag1.3
\endalign
$$
In general, any function that transforms like (1.2) (resp. (1.3))
is called a \it first-order \rm automorphic form (resp. \it
second-order \rm automorphic form).  Both Eisenstein series are
eigenfunctions of the hyperbolic Laplacian
$$
\Delta=-y^2\bigl(\frac{\partial^2}{\partial
x^2}+\frac{\partial^2}{\partial y^2}\bigr)=-4y^2 \frac{d
}{dz}\frac{d}{d \overline{z}},
$$
meaning
$$
\align
\Delta E(z,s)&=s(1-s)E(z,s),\\
\Delta E^*(z,s)&=s(1-s)E^*(z,s).
\endalign
$$
The second-order Eisenstein series $E^{*}(z,s)$ is known to have a
meromorphic continuation to all $s \in \C$ (see \cite{Gd},
\cite{O'S1}, \cite{Pe}). In \cite{G-O'S} it is shown that
$E^*(z,s)$ has a simple pole at $s=1$ with residue $-F(z)V^{-1}$,
meaning
$$
\lim\limits_{s \rightarrow 1} \left( E^*(z,s) +F(z)\frac{V^{-1}}{s-1}\right)
\,\,\,\,\,\text{\rm exists.}
$$
Recalling that the first-order Eisenstein series $E(z,s)$ has a
simple pole at $s=1$ with residue $V^{-1}$, we also can say that
$$
\lim\limits_{s \rightarrow 1} \left( E^*(z,s) +F(z)E(s,z)\right)
\,\,\,\,\,\text{\rm exists.}
$$
We take as the second-order analogue of the Kronecker limit
formula the study of the function
$$
\K_2(z)=\lim_{s\rightarrow 1}\left(E^*(z,s)+F(z)E(z,s)\right)
$$
for the following reason. By the definition of $F$, we have that
$$
E^*(z,s)+F(z)E(z,s)=\sum_{\g \in \G_\infty\backslash\G} F(\g
z)\Im(\g z)^s,
$$
which can be observed to be automorphic with respect to $\G$ for
all $s$, in particular when $s$ approaches $1$.  Therefore, the
function $\K_{2}$ is necessarily $\G$-invariant. Thus, in this notation,
$$
\lim\limits_{s \rightarrow 1} \left( E^*(z,s) +F(z)\frac{V^{-1}}{s-1}\right)
=\K_2(z)-F(z)\K_1(z).
$$

\vskip .10in Before describing our results concerning the
second-order Kronecker limit function $\K_{2}(z)$, we need the
following theorem concerning the first-order Kronecker limit
function $\K_{1}(z)$.

\proclaim{Theorem 1.1} The first-order Kronecker limit function
$\K_{1}$ admits the Fourier expansion 
$$
\K_1(z)=\sum_{n<0}k(n)e(n\overline{z}) + y+K-V^{-1}\log y +
\sum_{n>0}k(n)e(nz)
$$
with constants $K$ and $k(n)$. Furthermore,
$k(-n)=\overline{k(n)}$ and $k(n) \ll |n|^{1+\epsilon}$, with an
implied constant which depends solely on $\G$ and $\epsilon>0$.
\endproclaim

\vskip .10in We now can state the main results we obtain in our
study of the second-order Kronecker limit function $\K_{2}$.  To
begin, we have the following theorem regarding the convolution
Dirichlet series referred to in the title of the article.

\vskip .10in \proclaim{Theorem 1.2} Fix a positive integer $m$,
and let $k(0) = K + (\gamma + \log 4\pi m)/V$ where $K$ refers to a
component of the constant term in the Fourier expansion of
$\K_{1}$ and $V$ is the hyperbolic volume of $\G \backslash \H$.
Formally, for $s \in \C$, define the convolution Dirichlet series
$$
L^+_m(s)=\sum_{n=1}^\infty \frac{a_n k(m-n)}{n^s}
$$
and
$$
L^-_m(s)=\sum_{n=1}^\infty \frac{a_n}{n}\frac{
k(-m-n)}{(m+n)^{s-1}},
$$
which are formed from the Fourier coefficients of $f$ and $\K_{1}$.
Then the series $L^{+}_{m}$ and $L^{-}_{m}$ converge for $\Re (s)>
3$, admit a meromorphic continuation to all $s$ in $\C$ with $\Re (s)>1/2$, and are
holomorphic at $s=1$.
\endproclaim

\vskip .10in The usefulness of Theorem 1.2 will be evident in the
results below regarding the Fourier expansion of $\K_{2}$.

\vskip .10in It is known, and indeed is a elementary exercise,
that by combining the differential equation for $E(z,s)$ with its
Laurent expansion at $s=1$, one can prove the differential
equation $\Delta \K_{1}(z) = -V^{-1}$.  As we will see below, the
second-order analogue of this formula is the equation
$$
\Delta \K_2(z)=-8\pi i y^2 f(z) \frac{d}{d\overline{z}}\K_1(z).
$$
A more basic result would be to compute the differential equation
satisfied by $\frac d{dz}\K_{2}(z)$ or by $\frac d{d\overline{z}}\K_{2}(z)$.  We
carry out these derivations, ultimately proving the following two
theorems.

\vskip .10in \proclaim{Theorem 1.3} Let $\Pi_{hol}$ denote the
holomorphic projection operator for the space of smooth, weight
two automorphic forms into $S_2(\G)$.  Then
$$
\frac{1}{2\pi
i}\frac{d}{dz}\K_2(z)=f(z)\K_1(z)-\Pi_{hol}\bigl(f(z)\K_1(z)\bigr),
$$
Furthermore, if we set
$$
\K^+_1(z) =\sum\limits_{n>0}k(n)e(nz),
$$
then we have
$$
\Pi_{hol}\bigl(f(z)\K_1(z)\bigr)= \sum_{m=1}^\infty m L^+_m(1)
e(mz)-\frac{1}{2\pi i}F(z)\frac{d}{dz}\K_1^+(z)
+\frac{1}{4\pi}F(z).
$$
\endproclaim

\vskip .10in \proclaim{Theorem 1.4} Let $W_{s}(z)$ be the
classical Whittacker function associated to $\PSL_2(\R)$ and set
$$
W^{*}(z) = \frac{d}{ds}W_{s}(z)\Big|_{s=1} = \Gamma(0,4\pi
y)e^{4\pi y} e(z),
$$
where $\Gamma(s,a)$ denotes the incomplete gamma function
$$
\Gamma(s,a) = \int\limits_{a}^{\infty}e^{-t}t^{s-1}dt.
$$
Then
$$
\frac{d}{d\overline{z}}\K_2(z)=\frac{-2\pi i}V \sum_{n=1}^\infty
a_n W^*(nz)
+\frac{i}{2Vy}F(z)+F(z)\frac{d}{d\overline{z}}\K_1(z)+2\pi
i\sum_{m=1}^\infty mL^-_m(1)e(-m \overline{z}),
$$
\endproclaim

\vskip .10in As one would hope, Theorem 1.3 and Theorem 1.4
express the derivatives of $\K_{2}$ in terms of the initial information, namely
$\K_{1}$ and $f$.  Observe that either Theorem 1.3 or Theorem 1.4
can be used to compute $\Delta \K_{2}$; however, neither result
can be used to derive the other.

Theorem 1.3 is appealing because of its relatively concise
statement.  Theorems 1.3 and 1.4 indicate the necessity in studying the
Dirichlet series which are defined in Theorem 1.2.  At this point,
it must be noted that, in order to make sense out of Theorems 1.3 and 
1.4, we need to have some idea as to the growth of the special
values $L_{m}^{+}(1)$ and $L_{m}^{-}(1)$.  Before doing so, we state the following
result, which gives the complete Fourier expansion of the
second-order Kronecker limit function $\K_{2}$.

\vskip .10in \proclaim{Theorem 1.5} With notation as described
above, the second-order Kronecker limit function $\K_2(z)$ admits
the Fourier expansion
$$
\K_2(z)=\frac{-1}{V}\sum_{n=1}^\infty\frac{a_n}{n}W^*(nz)
-\sum_{m=1}^\infty L^+_m(1) e(mz)-\sum_{m=1}^\infty L^-_m(1)
e(-m\overline{z})+F(z)\K_1(z).
$$
\endproclaim

\vskip .10in Theorem 1.5 gives a complete description of the
second-order Kronecker limit function associated to
$E^{\ast}(z,s)$ at $s=1$.  The new ingredients that are not fully
understood are the special values $L_{m}^{+}(1)$ and
$L_{m}^{-}(1)$.  Theorem 1.2 asserts that $L^{+}_{m}(1)$ and
$L^{-}_{m}(1)$ are finite for all $m$, but to show that
 the Fourier expansion in Theorem 1.5 makes sense we bound the special values $L_{m}^{+}(1)$
and $L_{m}^{-}(1)$.  These bounds will imply that the series
expansions in Theorems 1.3, 1.4 and 1.5 converge for all $z \in
\H$.

\proclaim{Theorem 1.6} With the notation as above, we have the
bounds
$$
L^+_m(1), L^-_m(1) \ll  m^{1+\epsilon}
$$
with an implied constant that depends solely on $\G$, $f$ and
$\epsilon>0$. In addition, if: (i) the Fourier coefficients
$a_{n}$ of $f$ are in $\R$ for all $n$, and (ii) we have that
$\iota(\G)=\G$ where
$$
\left(\smallmatrix a & b \\ c & d \endsmallmatrix\right) \overset
\iota \to \longrightarrow \left(\smallmatrix -a & b \\ c & -d
\endsmallmatrix\right),
$$
then the special values $L^+_m(1)$ and $L^-_m(1)$ are also in $\R$
for all $m \geqslant 1$.
\endproclaim

To summarize, Theorem 1.1 establishes the Fourier expansion of the
first-order Kronecker limit function $\K_{1}$ and sets notation to
be used later.  Theorem 1.2 defines two families of convolution
Dirichlet series and asserts their meromorphic continuation and
holomorphicity at $s=1$.  Theorem 1.3 and Theorem 1.4 state two
different first-order differential equations which are satisfied
by the second-order Kronecker limit function $\K_{2}$, and Theorem
1.5 gives its Fourier expansion.  Bounds for the Fourier
coefficients of $\K_{1}$ are given in Theorem 1.1, and Theorem 1.6
gives analogous bounds for the Fourier coefficients of $\K_{2}$. We
believe that these results provide a complete investigation into
analytic aspects of the Fourier series development for $\K_{2}$.

\vskip .10in The outline of the paper is as follows.  In section 2
we initiate the development of the Fourier expansion of $\K_{2}$ and quickly find that
$$
\K_2(z)=A(z)+B(z)+F(z)\K_1(z)
$$
where $F\K_1$ is understood and $A$ is very similar to $F/V$ but non-holomorphic (in fact $\Delta A=F/V$). The main work in this paper is in understanding the term
$$
B(z)=\sum_{n=1}^\infty\bigl(b_ne(nz)+b_{-n}e(-n\bar{z})\bigr).
$$
The barrier to explicitly finding the constants $b_n$ is that they come from the Fourier coefficients $\phi^*_n(s)$ of the second-order Eisenstein series $E^*(z,s)$. These coefficients involve Kloosterman sums twisted by modular symbols and their values are not known inside the critical strip $0 \leqslant \Re (s) \leqslant 1$ even for the simplest congruence groups.

In section 3 we state, but do not prove, three key results: two on the analytic aspects of
Poincar\'e series, both holomorphic of weight 2 and non-holomorphic, and a third concerning
the holomorphic projection operator. Taken together these tools are powerful enough to probe the elements $b_n$. Because the proofs are so involved, we postpone verifying the statements of these
results 
until later in the paper.  

In section 4, we obtain information about the holomorphic part of $B$ by considering the holomorphic projection of the smooth, weight 2 function $\frac d{dz} \K_2(z)$. In the next section we show that the  coefficients $b_m$, for $m>0$, are given by the values of the convolution Dirichlet series $L^+_m(s)$ at $s=1$. This proves Theorem 1.3. A similar idea is used in section 6 to find the anti-holomorphic part of $B$ in terms of $L^-_m(1)$, proving Theorem 1.4. Combining these two theorems produces Theorem 1.5. There seems to be no symmetry between the holomorphic and anti-holomorphic parts of $B$. This is to be expected since the definition of $E^*$ includes a holomorphic cusp form $f$, breaking the symmetry.

In section 7 we complete the proof of Theorem 1.1
(bounding the Fourier coefficients of $\K_{1}$)
and prove Theorem 1.6 concerning the bounds on $b_n$.  All results in section 7 come
from careful considerations involving the asyptotics of $E(z,s)$ and $E^*(z,s)$ as $z$ approaches cusps. The crude bounds coming from the meromorphic continuation of these series are improved by a type of bootstrapping procedure. These results are independent of those
 in previous sections.  At this time, there are a few
remaining pieces to complete:  The proofs of the results in section 3 as well as the meromorphic continuations
and regularity at $s=1$ of $L^{+}_{m}(s)$ and $L_{m}^{-}(s)$.  In
section 8 we use the spectral theory of automorphic forms to prove Theorem 3.1 and Theorem 3.2, and in section 9
we prove Proposition 3.3 as well as the remaining properties
regarding $L^{+}_{m}$ and $L^{-}_{m}$ by introducing a type of non-holomorphic Poincar\'e series, $Q_m(z,s;F)$, that includes $F$ in its definition.  

Finally, in section 10 we
conclude with two types of examples:  The first example shows how
to explicitly evaluate the first-order Kronecker limit function
$\K_{1}$ for the congruence subgroups $\G_{0}(N)$ for square-free
$N$, and the second example poses, as well as numerically
investigates, a problem related to Theorem 1.3 involving the holomorphic projection
operator.

\vskip .10in The detailed, technical results in this paper begin
in section 7, then carry through to sections 8 and 9. These
precise calculations are used to prove the statements in section 3 and the
meromorphic continuation of $L_{m}^{+}$ and $L_{m}^{-}$, the details of which
comprise the most difficult parts of our work.
  The arrangement of sections in
this paper is meant to provide the motivation for each new result as it is needed and is purposefully consistent with our order of discovery.

\vskip .25in \head \bf \S 2. The Fourier expansion of $\K_{2}$
\endhead

Our starting point is the Fourier expansions for the functions
$E^{\ast}(z,s)$, $E(z,s)$ and $F(z)$, from which we obtain a
somewhat general Fourier expansion for $\K_{2}(z)$.  From
\cite{O'S1}, page 164, we have that the Fourier expansion for the
second-order Eisenstein series $E^{\ast}(z,s)$ is
$$
E^*(z,s)=\sum_{n\neq 0}\phi^*_n(s)W_s(nz)\tag 2.1
$$
where $W_s$ is the Whittaker function
$$
W_{s}(nz) = 2 \vert n \vert ^{1/2}y^{1/2}K_{s-1/2}(2\pi \vert n
\vert y) e( n x)
$$
and $K_{s}$ is the $K$-Bessel function
$$
K_{s}(z) = \frac{1}{2}\int\limits_{0}^{\infty}e^{-z(u +
1/u)/2}u^{s} \frac{du}{u} \,\,\,\,\,\text{\rm for $\Re (s) > 0$.}
$$
Note that we have also used Corollary 4.3 of \cite{O'S1} which
proves that, in this instance, the second-order Eisenstein series
has no constant term in its Fourier expansion. Exact formulas in
terms of number theoretic functions are known for the Fourier
coefficients of the first-order Eisenstein series $E(z,s)$ in the
case when $\G$ is a congruence subgroup.  In general, no such
formulas are known for the coefficients $\phi^*_n(s)$.  Let us use
the following Laurent series:
$$
\align
E(z,s)&=\frac{V^{-1}}{s-1}+\K_1(z)+O(s-1),\\
\phi^*_n(s)&=\frac{b_n(-1)}{s-1}+b_n(0)+O(s-1),\\
W_s(z)&=e(z)+W^*(z)(s-1)+O((s-1)^2)
\endalign
$$
where $W^*(z)=\left.\frac{d}{ds}W_s(z)\right|_{s=1}$.  In
Corollary 2.2 below, we will prove the formula for $W^*(z)$
asserted in Theorem 1.4. Note that, by definition,
$W_s(z)=W_s(\bar{z})$ for $z$ in the lower half plane.
Substituting these expansions into the defintion
$$
\K_2(z)=\lim_{s\rightarrow 1}\left(E^*(z,s)+F(z)E(z,s)\right)
$$
yields the expression
$$
\multline
\K_2(z)=\lim_{s\rightarrow 1}\left[\sum_{n\neq 0}\frac{b_n(-1)}{s-1}e(nz)+F(z)
\frac{V^{-1}}{s-1}\right.\\
\left.+\sum_{n=1}^\infty\bigl(b_n(0)e(nz)+b_{-n}(0)e(-n\bar{z})\bigr)
+\sum_{n\neq 0}b_n(-1)W^*(nz)+F(z)\K_1(z)\right].
\endmultline
$$
Since the limit which defines $\K_{2}(z)$ exists, it is evident
that we must have
$$
b_n(-1)= \left\{ \matrix \displaystyle \frac{-a_n}{n}V^{-1} & n
\geqslant 1 \\  \\  \displaystyle 0 & \text{\rm otherwise} \endmatrix
\right\}.
$$
Set $b_{n} = b_{n}(0)$ and, at this time, we can write
$$
\K_2(z)=A(z)+B(z)+F(z)\K_1(z)\tag 2.2
$$
where
$$
\align
A(z)&=\frac{-1}{V}\sum_{n=1}^\infty\frac{a_n}{n}W^*(nz),\tag 2.3\\
B(z)&=\sum_{n=1}^\infty\bigl(b_ne(nz)+b_{-n}e(-n\bar{z})\bigr).\tag 2.4
\endalign
$$
To go further, we compute the resulting formula obtained by
applying the Laplacian to $\K_{1}$, $\K_{2}$, $A$, and $B$. Using
that
$$
\Delta\left[\lim_{s\rightarrow1}
\left(E(z,s)-\frac{V^{-1}}{s-1}\right)\right]= \lim_{s\rightarrow
1}\left[\Delta\left(E(z,s)-\frac{V^{-1}}{s-1}\right)\right] =
\lim_{s\rightarrow
1}\left[s(1-s)E(z,s)\right]
$$
one shows that
$$
\Delta \K_1(z)=-V^{-1}.\tag 2.5
$$
Similarly, we now consider
$$
\Delta\lim_{s\rightarrow 1}\left(E^*(z,s)+F(z)E(z,s)\right)=\lim_{s\rightarrow 1}
\Delta\left(E^*(z,s)+F(z)E(z,s)\right),
$$
which can be easily computed.  Since
$$
\Delta \left(E^*(z,s)+F(z)E(z,s)\right) =
s(1-s)\left(E^{\ast}(z,s) + F(z)E(z,s)\right) - 8 \pi iy^{2}f(z)
\frac{d}{d \overline{z}}E(z,s),
$$
we then obtain, by taking $s \rightarrow 1$, the formula
$$
\Delta \K_2(z) =-8\pi i y^2 f(z) \frac{d}{d\overline{z}}\K_1(z).
$$
Also $\Delta B(z)=0$, so then we have by $(2.2)$ that
$$
\Delta \K_2(z) = \Delta A(z) + \Delta ( F(z)\K_{1}(z)),
$$
which, when combined with the above formulas, yields the
expression
$$
\Delta A(z) =F(z)V^{-1}.\tag 2.6
$$

\vskip .10in In order to examine $A(z)$ more explicitly, we shall
study $W^*(nz)$ by means of its definition in terms of
$K$-Bessel functions. For this, we use that the $K$-Bessel
function can be written as
$$
K_{s-1/2}(2\pi y)=\frac{\sqrt{\pi}}{\G(s)}(\pi y)^{s-1/2}\int_1^\infty (t^2-1)^{s-1}
e^{-2\pi ty}\,dt.\tag 2.7
$$
(see page 205, \cite{Iw1}).  The integral in (2.7) converges
absolutely for $\Re (s) >0$ and $y>0$. We want to find
$\displaystyle W^*(z)=\left.\frac{d}{ds}W_s(z)\right|_{s=1}$.

\vskip .10in
\proclaim{Lemma 2.1}  For all $y > 0$, we have
$$
\left.\frac{d}{ds}K_{s-1/2}(2\pi y)\right|_{s=1}=
\G(0,4\pi y)\frac{e^{2\pi y}}{2\sqrt{y}}.
$$
\p Trivially, we have
$$
\align \displaystyle \left.\frac{d}{ds}\int_1^\infty
(t^2-1)^{s-1}e^{-2\pi ty}\, dt\right|_{s=1} & \displaystyle
=\int_1^\infty e^{-2\pi ty}\log(t-1)\, dt+\int_1^\infty e^{-2\pi
ty}\log(t+1)\, dt\\  \\ &\displaystyle =e^{-2\pi y}\int_0^\infty
e^{-2\pi uy}\log u\, du+e^{2\pi y}\int_2^\infty e^{-2\pi uy}\log
u\, du.
\endalign
$$
Now
$$
\int_0^\infty  e^{-u}\log u\, du=\G'(1),
$$
where $\Gamma(s)$ denotes the
classical gamma function.  Therefore,
$$
\int_0^\infty e^{-2\pi uy}\log u\, du=\frac{1}{2\pi y}(\G'(1)-\log 2\pi y).
$$
Through elementary computations, using integration by parts, one
can show that
$$
\int_2^\infty e^{-2\pi uy}\log u\, du=\frac{1}{2\pi y}(\G(0,4\pi
y)+ e^{-4\pi y}\log 2).
$$
Combining these formulas, we obtain the relation
$$
\left.\frac{d}{ds}\int_1^\infty (t^2-1)^{s-1}e^{-2\pi ty}\,
dt\right|_{s=1} =\frac{1}{2\pi y}\left(\G(0,4\pi y)e^{2\pi
y}+(\G'(1)-\log \pi y)e^{-2\pi y}\right).
$$
To complete the proof, one simply computes the derivative of (2.7)
with respect to $s$ and sets $s=1$.  Using that
$$
\int\limits_{1}^{\infty}(t^{2}-1)^{s-1}e^{-2\pi t y}dt \Big|_{s=1}
= \frac{e^{2\pi y}}{2 \pi y},
$$
the result follows from the standard rules of calculus.
\bs

\vskip .10in \proclaim{Corollary 2.2} For $z$ in $\H$ and $n \geqslant
1$, we have the following formulas:
$$
\align
W^*(nz)=& \G(0,4\pi ny) e^{4\pi ny} e(nz),\\
\frac{d }{d z}W^*(nz)=& \frac{i}{2y}e(nz),\\
\frac{d }{d \overline{z}}W^*(nz)=&
\frac{-i}{2y}e(nz)+2\pi inW^*(nz), \\
\Delta W^{\ast}(nz) =& - e(nz).
\endalign
$$
\p The first identity follows directly from the definition of the
Whittacker function in terms of the $K$-Bessel function, together
with Lemma 2.1 and elementary calculus.  The remaining three
formulas are direct computations from the first expression, using
nothing more than the fundamental theorem of calculus and standard
formulas for differentiation of functions of one complex variable.
\bs

These computations allow us to give a precise description of
$A(z)$. Indeed, by definition we have
$$
A(z)=\frac{-1}{V}\sum_{n=1}^\infty\frac{a_n}{n}W^*(nz),
$$
so then Corollary 2.2 allows one to compute various derivatives of
$A(z)$.

\vskip .25in \head \bf \S 3. Poincar\'e series and holomorphic
projection \endhead

In order to continue studying the computations given in the
previous section, we will use the holomorphic projection operator,
whose basic properties we recall in the present section.

For any two smooth functions $\varphi_1,\varphi_2$ which transform
with weight $k$, and have exponential decay at the cusps, the
Petersson inner product between $\varphi_{1}$ and $\varphi_{2}$ is
defined by
$$
\langle \varphi_1,\varphi_2 \rangle_k =\int_{\G \backslash \H} y^k
\varphi_1(z)\overline{\varphi_2(z)}\,d\mu (z),
$$
where $d\mu (z)=dx dy/y^2$ is the usual hyperbolic volume form. It
can be shown that the Petersson inner product is non-degenerate on
the space of holomorphic weight $k$ cuspforms $S_k(\G)$.
Consequently, for any $\varphi_1$ as above, there exists a form
$\Pi_{hol}(\varphi_1)$ in $S_k(\G)$ such that for every $g$ in
$S_k(\G)$
$$
\langle \varphi_1,g \rangle_k=\langle
\Pi_{hol}(\varphi_1),g\rangle_k.
$$
The image $\Pi_{hol}(\varphi_1)$ of $\varphi_{1}$ into $S_{k}(\G)$
is called the \it holomorphic projection \rm of $\varphi_1$.

In the appendix of \cite{Za1}, beginning on page 286, it is shown
that the Fourier coefficients of $\Pi_{hol}(\varphi_1)$ can be
computed by taking $g$ in the inner product above to be the weight
$k$ holomorphic Poincar\'e series. In section 8 we construct and
study aspects of these series relevant for our work.  In order to
proceed with the computations from the previous section, we shall
state here various results regarding these Poincar\'e series,
leaving the proofs of the assertions until section 8.

Thus far we have only concerned ourselves with a single cusp,
which we assumed was uniformized to be at the point at $\infty$.
Let us now consider the possibility that an arbitrary (finite)
number of $\G$-inequivalent cusps exists.  If there are other
inequivalent cusps,  let us fix representatives, label them
$\ca,\cb,\cc \dots$ and use the scaling matrices $\sa,\sb,\sc
\dots$ to give local coordinates near these cusps (see Chapter 2
of \cite{Iw1} as well as \cite{O'S1}). The subgroup $\G_\ca$ is the
set of elements of $\G$ which fixes the cusps equivalent to $\ca$,
and we have that
$$
\sa^{-1} \G_\ca \sa= \G_\ci=
\left\{ \pm \left(\smallmatrix 1 & m \\ 0 & 1
\endsmallmatrix\right)
\; \big | \; \ m\in {\Bbb Z}\right\}.
$$

\noindent Following Selberg \cite{Se}, for each $m \geqslant 1$, we
define the non-holomorphic Poincar\'e series associated to the
cusp $\ca$ as
$$
U_{\ca m}(z,s)=\sum_{\g \in \G_\ca\backslash\G} \Im(\sa^{-1}\g
z)^s e(m\sa^{-1}\g z).
$$
We shall also need $U'_{\ca m}=\frac{d}{dz}U_{\ca m}(z,s)$, the termwise derivative of $U_{\ca m}$.

\vskip .10in To examine the growth of $U_{\ca m}$ and other automorphic functions, we follow the
convention set in (2.42) of \cite{Iw1} and introduce the
useful notation
$$
y_\G(z)=\max_\ca \max_{\g \in \G}(\Im( \sa^{-1}\g z)).
$$
Heuristically, the function $y_{\G}(z)$ measures how close the
point $z\in \G\backslash \H$ is to a cusp. If $\psi$ (or $|\psi|$)
is a smooth weight zero form (i.e., $\G$-invariant function),
then it is more convenient to write
$$
\psi(z) \ll y_\G(z)^A,
$$
than, for example, writing that $\psi(\sa z) \ll y^A$ for each
cusp $\ca$ as $y \rightarrow \infty$.

\vskip .10in \proclaim{Theorem 3.1} For all $m \geqslant 1$ and $\Re (s)
> 1$, the series $U_{\ca m}(z,s)$ and $\frac{d}{dz}U_{\ca m}(z,s)$ are
pointwise absolutely convergent  and uniformly convergent for $s$ in
 compact sets.  Furthermore, both series admit
meromorphic continuations to all $s \in \C$ which are analytic at
$s=1$. For $\Re (s)> 1/2$ we have the growth conditions
$$
U_{\ca m}( z,s) \ll |m|^{-1/2}\sqrt{y_\G(z)}
$$
and
$$
yU'_{\ca m}( z,s) \ll  |m|^{-1/2}\sqrt{y_\G(z)}
$$
with an implied constant depending on $s$ and $\G$ alone.
\endproclaim

\vskip .10in Going further, let us define
$$
V_{\ca m}(z,s)=\sum_{\g \in \G_\ca\backslash\G} \Im(\sa^{-1}\g
z)^s e(m\sa^{-1}\g z) j(\sa^{-1}\g, z)^{-2}
$$
which can be viewed as a weight two version of $U_{\ca m}$.
Formally, we would like to define our weight two holomorphic
Poincar\'e series, which we will denote by $P_{\ca m}(z)_2$, to be
given by $V_{\ca m}(z,0)$.  However, as will be evident from the
analysis of section 8, the series defining $V_{\ca m}(z,s)$ is
absolutely convergent only for $\Re (s)>0$.  In order to address
this difficulty, we proceed as follows.

\vskip .10in By a direct computation, one can easily show that for
any $z \in \H$ and $\gamma \in \PSL_2(\R)$, we have the identity
$$
2i\frac{d}{dz} \left[\Im (\g z)^s e(m\g z)\right]=s\frac{\Im (\g
z)^{s-1}}{j(\g,z)^2}e(m\g z) -4\pi m \frac{\Im (\g
z)^{s}}{j(\g,z)^2}e(m\g z).
$$
By summing over all coset representatives $\gamma \in \G_\ca\backslash\G$, this implies the
formula
$$
sV_{\ca m}(z,s-1)=2i\frac{d}{dz}U_{\ca m}(z,s)+4\pi mV_{\ca
m}(z,s),\tag 3.1
$$
which necessarily holds only in the half-plane of absolute
convergence for both series which define $U_{\ca m}$ and $V_{\ca
m}$.  Therefore, in the light of Theorem 3.1, it makes (formal)
sense to define the Poincar\'e series $P_{\ca
m}(z)_2$ through the formula
$$
P_{\ca m}(z)_2=2i\frac{d}{dz} U_{\ca m}(z,1) +4\pi m V_{\ca m}(z,1).
$$
We verify in Theorem 3.2 below that this does indeed give us a weight two holomorphic cusp form.
Let us now examine how one can evaluate various inner products involving $P_{\ca m}(z)_2$.

\vskip .10in Given a suitable function $\varphi$, we propose to
evaluate $\langle \varphi ,P_{\ca m}(\cdot)_2 \rangle_2$ by first
studying the meromorphic function
$$
\langle \varphi , V_{\ca m}(\cdot,\overline{s}-1) \rangle_2
=\int_0^\infty \int_0^1 \varphi(z) y^{s-1} \overline{e(mz)} \,dx\,dy. \tag 3.2
$$
Under certain restrictions on $\varphi$ the unfolded inner product on the right of (3.2) will converge for
$\Re (s) $ large and may be computed to yield a function with a natural meromorphic continuation (for example involving gamma functions) to $s=1$. In this way (3.2) at $s=1$ yields an
evaluation of $\langle \varphi ,P_{\ca m}(\cdot)_2 \rangle_2$.
Indeed, we will follow this method to prove the
following theorem.

\vskip .10in \proclaim{Theorem 3.2} The weight two  Poincar\'e
series $P_{\ca m}(z)_2$ is in $S_2(\G)$, the vector space of
holomorphic weight two cusp forms with respect to $\G$.
Furthermore, for any $f$ in $S_2(\G)$ with
$$
j(\sa, z)^{-2}f(\sa z)=\sum_{n=1}^\infty a_\ca (n) e(nz)
$$
we have that
$$
\langle f,P_{\ca m}(\cdot)_2 \rangle_2=a_\ca (m)/(4\pi m).
$$
\endproclaim

\vskip .10in Frequently, we will assume that the cusp in question
has been uniformized to be at $\ci$, so then, for ease of
notation, we will set $U_m=U_{\ci m}$, $V_m=V_{\ci m}$ and
$P_m=P_{\ci m}$.  From the above discussion, we have the
following.  If $\varphi$ is a smooth, bounded, continuous function
on $\H$ which transforms like a weight two form with respect to
the action by $\G$, we then have
$$
\Pi_{hol}(\varphi)=\sum_{m=1}^\infty d_m e(mz) \text{ \ \ where \
\ } d_m=4\pi m \langle \varphi,P_{m}(\cdot)_2 \rangle_2.\tag 3.3
$$

\noindent
 As stated above, the proofs of Theorem 3.1 and Theorem
3.2 will be given in section 8 below.

\vskip .10in In the forthcoming work, we will make use of the
following proposition.

\vskip .10in \proclaim{Proposition 3.3} Let $\varphi_{1}$ be a
smooth weight zero form (function) and $\varphi_{2}$ a smooth form
of weight two.  Then:

\roster
\item"{(i)}"
The form $\frac{d}{dz}\varphi_{1}$ is a weight two form;
\item"{(ii)}"
The form $y^{2}\frac{d}{d\overline{z}}\varphi_{2}$ is a weight
zero form;
\item"{(iii)}"
Assuming appropriate growth conditions on the functions near the
cusps, we have the inner product formula
$$
\langle \frac{d}{dz}\varphi_1, \varphi_2 \rangle_2=-\langle
\varphi_1, y^2\frac{d}{d\overline{z}}\varphi_2 \rangle_0.
$$
\endroster
The growth conditions are satisfied, for example, if $\varphi_1$
and $\frac{d}{dz}\varphi_1$ have at
most polynomial growth in $y$ in the cusps and if $\varphi_2$ and
$y^2\frac{d}{d\overline{z}}\varphi_2$ have exponential decay in the cusps.
\endproclaim

\vskip .10in Proposition 3.3 will be proved as a corollary to
Proposition 9.3, which states a more general result involving the
Maass weight raising and lowering operators.  We state the
specific result here in order to continue with the calculations
given in section 2.  We note that the proof of Proposition 2.1.3
of \cite{Bu}, which involves Stokes's theorem, may be adapted to
yield a proof of Proposition 3.3.  Rather than following this
approach, our proof of Proposition 9.3 involves integration by
parts together with some aspects of the first-order Eisenstein
series, which gives an argument that extends to consider others
pairs of forms with complementary weights.

\vskip .10in Directly from Proposition 3.3, we have the following.

\vskip .10in \proclaim{Corollary 3.4} Let $\varphi$ be a
smooth, weight $0$ function on $\H$ which is $\G$ invariant. Assume
that $\varphi$ and $\frac{d}{dz}\varphi$ have at most
polynomial growth in the cusps of $\G \backslash \H$.  Then
$$
\Pi_{hol}\left(\frac{d}{dz}\varphi\right)=0.
$$
\p From Theorem 3.2, we have that the weight two Poincar\'e series
is holomorphic, i.e.
$$
\frac{d}{d\overline{z}} P_m(z)_2=0.
$$
Corollary 3.4 now follows by using the second part of Theorem
3.2 together with Proposition 3.3. \bs

\vskip .10in To re-iterate, the proofs of Theorem 3.1 and Theorem
3.2 will be given in section 8, and the proof of Proposition 3.3
will be given in section 9.

\vskip .25in \head \bf \S 4. $\K_{2}$ and the holomorphic
projection of $f\K_{1}$ \rm
\endhead

\vskip .10in Using the material stated in section 3, we now
continue with the calculations from section 2.  Specifically, we
will complete the proof of Theorem 1.3 in this section and the next.

\vskip .10in Recall from (2.2) that we have written
$$
\K_{2}(z) = A(z) + B(z) + F(z)\K_{1}(z),
$$
with $A(z)$ and $B(z)$ defined in (2.3) and (2.4),
respectively.  Using Corollary 2.2, we then get the formula
$$
\align
\frac{d}{dz}\K_2(z)&=\frac{d}{dz}A(z)+\frac{d}{dz}B(z)+\frac{d}{dz}(F(z)\K_1(z))\\
&=\frac{-i}{2Vy}F(z)+2\pi i\sum_{n=1}^\infty n b_n  e(nz)+
F(z)\frac{d}{dz}\K_1(z)+2\pi i f(z)\K_1(z).\tag4.1
\endalign
$$
The right-hand side of (4.1) is a sum of two weight two forms. Since
the holomorphic projection operator is linear, we then have
$$
\align
\Pi_{hol}\left(\frac{d}{dz}\K_2(z)\right)&=\Pi_{hol}\left(\frac{-i}{2Vy}F(z)+
2\pi i\sum_{n=1}^\infty n b_n  e(nz)
+F(z)\frac{d}{dz}\K_1(z)\right)
\\&+\Pi_{hol}\left(2\pi i f(z)\K_1(z)\right).
\endalign
$$
By Corollary 3.4, as will be established with the verification of
the appropriate growth conditions in the proof given of
Proposition 9.3, we have
$$
\Pi_{hol}\left(\frac{d}{dz}\K_2(z)\right) = 0
$$
so then
$$
\Pi_{hol}\left(\frac{-i}{2Vy}F(z)+ 2\pi i\sum_{n=1}^\infty n b_n
e(nz) +F(z)\frac{d}{dz}\K_1(z)\right)  +\Pi_{hol}\left(2\pi i
f(z)\K_1(z)\right) = 0 .
$$
Let
$$
g(z)=\frac{-i}{2Vy}F(z)+2\pi i\sum_{n=1}^\infty n b_n
e(nz)+F(z)\frac{d}{dz}\K_1(z). \tag 4.2
$$
We now will show that (4.2) is actually a holomorphic cusp form,
and hence equal to its own holomorphic projection.  Therefore
 $g=-2\pi i\, \Pi_{hol}(f\K_1)$ and substituting this back into (4.1) will complete the proof of
 the first part of Theorem 1.3.

\vskip .10in Using the differential equation $(2.5)$ for $\K_{1}$, we get
$$
\align \frac{d}{d \bar z} g(z) & = \frac{d}{d \bar
z}\left(\frac{-i}{2Vy}F(z) + F(z) \frac{d}{dz}\K_{1}(z)\right)
\\&= \frac{-i}{2V}F(z)\frac{d}{d \bar z}(y^{-1}) + F(z) \frac{d^{2}}{dz d \bar z}\K_{1}(z)
\\&= \frac{-i}{2V}F(z)\frac{1}{2i}y^{-2} + F(z)\left(-\frac{1}{4y^{2}}(-V^{-1})\right) = 0,
\endalign
$$
hence $g$ is holomorphic.  It thus remains to show that $g$ has
exponential decay in each cusp, which will follow by studying its Fourier
expansion in each cusp.  In this generality, there are a number of
analytic quantities associated with the cusp $\ca$.  Using an
obvious extension of notation established thus far, we define:
$$
\align
E_\ca(z,s)&=\sum_{\g \in \G_\ca\backslash\G} \Im(\sa^{-1}\g z)^s,\\
E^*_\ca(z,s)&=\sum_{\g \in \G_\ca\backslash\G} \langle \g, f\rangle \Im(\sa^{-1}\g z)^s,\\
F_\ca(z)&=2\pi i\int_\ca^z f(w) \, dw,\\
\K_{1\ca}(z)&=\lim_{s\rightarrow 1}\left(E_\ca(z,s)-\frac{V^{-1}}{s-1}\right),\\
\K_{2\ca}(z)&=\lim_{s\rightarrow 1}\left(E_\ca
^*(z,s)+F_\ca(z)E_\ca(z,s)\right).
\endalign
$$
The relevant Fourier expansions at the cusp $\cb$ are:
$$
\align
E_\ca(\sb z,s)&=\delta_{\ca \cb}y^s +\phi_{\ca \cb}(s)y^{1-s}+
\sum_{n\neq 0}\phi_{\ca \cb}(n,s)W_s(nz),\tag 4.3\\
E_\ca^*(\sb z,s)&=\phi^*_{\ca \cb}(0,s)y^{1-s}+\sum_{n\neq 0}
\phi^*_{\ca \cb}(n,s)W_s(nz),\tag 4.4\\
j(\sb,z)^{-2}f(\sb z)&=\sum_{n=1}^\infty a_\cb(n) e(nz),\\
F_\ca(\sb z)&= T_{\ca \cb}+ \sum_{n=1}^\infty \frac{a_\cb(n)}{n}
e(nz),\tag 4.5
\endalign
$$
where we define the period $T_{\ca \cb} = 2\pi i
\int_\ca^\cb f(w)\, dw$. We refer to equation (3.20) \cite{Iw1}
for a proof of (4.3), and to equation (1.1) of \cite{O'S1} for a
proof of (4.4). Note that by Corollary 4.3 of \cite{O'S1} we have
that $\phi^*_{\ca \ca}(0,s)=0$ which agrees with (2.1). We write
the Laurent expansion of $\phi^*_{\ca \cb}(n,s)$ at $s=1$ as
$$
\phi^*_{\ca \cb}(n,s)=\frac{b_{\ca \cb}(n,-1)}{s-1}+b_{\ca
\cb}(n,0)+O(s-1).
$$
The analogue of (2.2) for $\K_{2\ca}$ at the cusp $\cb$ is then
$$
\align \K_{2\ca}(\sb z)&=\frac{T_{\ca \cb}}{V} \log y+ b_{\ca
\cb}(0,0)
+\frac{-1}{V}\sum_{n=1}^\infty \frac{a_\cb(n)}{n}W^*(nz) \\
&+ \sum_{n=1}^{\infty} \bigl(b_{\ca \cb}(n,0)e(nz)+b_{\ca
\cb}(-n,0)e(-n\overline{z})\bigr) + F_\ca(\sb z) \K_{1\ca}(\sb z).
\endalign
$$
Assuming the Fourier expansion
$$
\K_{1\ca}(\sb z)=\sum_{n<0}k_{\ca \cb}(n)e(n\overline{z}) +
\delta_{\ca \cb}y+k_{\ca \cb}(0) -V^{-1}\log y +\sum_{n>0}k_{\ca
\cb}(n)e(nz),\tag4.6
$$
which we will establish in this section below, we see that
$$
\align
\frac{d}{dz} \K_{2\ca}(\sb z)& = 2\pi i\sum_{n=1}^\infty n b_{\ca \cb}(n,0)e(nz) \\
&+  F_\ca (\sb z) \left(\frac{-i}{2} \delta_{\ca \cb} +
\frac{d}{dz}\left(\sum_{n>0}k_{\ca \cb}(n)e(nz)\right)\right)
+j(\sb ,z)^{-2}f(\sb z)\K_{1\ca}(\sb z).
\endalign
$$
Therefore, it follows that
$$
\align j(\sb, z)^{-2}&g(\sb z)=\frac{d}{dz}\K_2(\sb z)-2\pi i
\cdot j(\sb, z)^{-2}f(\sb z)\K_1(\sb z)\\& =2\pi
i\sum_{n=1}^\infty n b_{\ca \cb}(n,0)e(nz) +  F_\ca (\sb z)
\left(\frac{-i}{2} \delta_{\ca \cb} +
\frac{d}{dz}\left(\sum_{n>0}k_{\ca \cb}(n)e(nz)\right)\right).
\tag 4.7
\endalign
$$
If $\ca \neq \cb$, then (4.7) has rapid decay, which is seen
by combining (4.5) together with the fact that $\delta_{\ca \cb} =
0$.  If $\ca = \cb$, then in (4.5) we have that $T_{\ca \cb} =
0$, so we again conclude that (4.7) has rapid decay.

\vskip .10in Finally, it needs to be verified that the expansion
(4.6) holds.  Indeed, this expansion follows directly from the
Fourier expansion for the first-order non-holomorphic Eisenstein
series, as stated in (4.3), together with the special function
calculations given in the proof of Corollary 2.2.  The important
point is that the coefficient of $\log y$ is $V^{-1}$, which is
implied by fact that the Eisenstein series (4.3) has a first order
pole at $s=1$ with residue equal to $V^{-1}$.  All of these
properties of the Eisenstein series are proved in, for example,
\cite{Iw1} and \cite{Kub}.  This argument, which follows the
method of calculation given in section 2, gives the first part of
Theorem 1.1. The bounds on the Fourier coefficients of $\K_1$, claimed in the second part of Theorem 1.1, are achieved in section 7.

\vskip .10in With all this, the proof of the first statement of Theorem 1.3 is complete,
and, indeed, we have shown that for any cusp $\ca$, 
$$
\frac{1}{2\pi
i}\frac{d}{dz}\K_{2\ca}(z)=f(z)\K_{1\ca}(z)-
\Pi_{hol}\bigl(f(z)\K_{1\ca}(z)\bigr).
$$
It remains to give the stated expression for $\Pi_{hol}\bigl(f(z)\K_{1\ca}(z)\bigr)$ in the second part of Theorem 1.3. This is carried out next.

\vskip .25in \head \bf \S 5. The Dirichlet series $L_{m}^{+}$ and
the holomorphic projection of $f\K_{1}$
\endhead

\vskip .10in We continue by studying the Fourier coefficients of
$\Pi_{hol}\bigl(f(z)\K_{1\ca}(z)\bigr)$.  As in the previous section, we
will use the results stated in section 3, whose proofs we will
give in section 8.  In the notation established in section 1, let
us write
$$
L^{++}_m(s)=\sum_{n=m+1}^\infty \frac{a_n k(m-n)}{n^s},\tag 5.1
$$
so then, referring to the notation from Theorem 1.2, we have
$$
L^+_m(s)=\sum_{n=1}^{m-1} \frac{a_n
k(m-n)}{n^s}+\frac{a_m}{m^s}\left(K+\frac{\gamma+\log 4\pi
m}{V}\right)+L^{++}_m(s).
$$

\vskip .10in \proclaim{Proposition 5.1} With notation as above and
for $\Re (s) $ sufficiently large, we have the identity
$$
\align \langle f\K_1, V_m(\cdot,\overline{s}-1)\rangle_2 &
=\frac{\G(s)}{(4\pi)^s} L^{++}_m(s)+ \frac{\G(s)}{(4\pi
m)^s}\sum_{l=1}^{m-1} a_l k(m-l)\\& +\frac{a_m}{(4\pi m)^s}\left(
\frac{\G(s+1)}{4\pi m}+K \G(s) +\frac{-\G'(s)+\G(s)\log 4\pi m}{V}
\right).
\endalign
$$

\p This is carried out using the ideas of section 3, in particular
(3.2), and follows the line of standard computations.  First,
expand $f$ and $\K_{1}$ in their Fourier expansions, i.e.
$$
f(z) \K_{1}(z) = \left(\sum\limits_{n=1}^{\infty}a_{n}e(nz)\right)
\left(\sum_{n<0}k(n)e(n\overline{z}) + y+K-V^{-1}\log y +
\sum_{n>0}k(n)e(nz) \right).
$$
Next, unfold the integral in question, similar to (3.2), and carry
out the integral, ultimately using standard formulas for the
classical $\Gamma$ function.
\bs

\vskip .10in {\bf Remark 5.2.} The trivial bound for coefficients
of a weight two cusp form states that $a_n \ll n$ (see (5.7) of
\cite{Iw2}).  In section 7 below we will prove, as asserted in
Theorem 1.1, that $k(n) \ll |n|^{1+\epsilon}$. Therefore, it
follows that $L^+_m(s)$ is absolutely and uniformly convergent for
$\Re (s)>3$, as claimed in Theorem 1.2.  The meromorphic
continuation of $L^{+}_{m}(s)$ will follow from the expression
derived in Proposition 5.1 together with a study of the Poincar\'e
series $V_m(\cdot,s)$.

\vskip .10in We now work with the expression
$$
\langle f\K_1, P_m(\cdot)_2\rangle_2=\left.\langle f\K_1,
V_m(\cdot,\overline{s}-1)\rangle_2\right|_{s=1},
$$
in order to compute the Fourier coefficients of
$\Pi_{hol}(f(z)\K_{1}(z))$.  To begin, let us make sure that we
have, or will, establish enough results regarding $V_m(z,s)$ to
proceed. Assuming Theorem 3.1, from which we obtained (3.1), then
$V_m(z,s)$ has a meromorphic continuation to all $s \in \C$.
Again, Theorem 3.1 will be proved in section 8, at which time it
also will be shown that $V_m(z,s)$ has at most polynomial growth
in $y$ at the cusps. Therefore, $\langle f\K_1,
V_m(\cdot,\overline{s}-1)\rangle_2$ converges to a meromorphic
function. Proposition 5.1 holds for $\Re (s) > 3$, as stated in
Remark 5.2, so then we now have the meromorphic continuation of
the Dirichlet $L^{++}_m(s)$, and hence $L^{+}_m(s)$, to all $s \in
\C$. The continuation $L^{+}_m(s)$ will not have a pole at $s=1$
once it has been shown that $V_m(z,s-1)$ does not have a pole at
$s=1$.  Hence, the stated results in section 3, together with the
growth condition for $V_m(z,s)$ which also comes from section 8,
are sufficient to allow us to continue our calculations.

\vskip .10in  Direct calculations using the Fourier expansion of
$\K_{1}$ show that
$$
\frac{d}{dz} \K_{1}(z) = \frac{d}{dz} \K^{+}_{1}(z) +
\frac{i}{2Vy} - \frac{i}{2}.
$$
Combining this with equation (4.2), as well
as subsequent discussion, we get
$$
\align
\Pi_{hol}(f(z)\K_{1}(z))&=\frac{1}{2\pi i}\left(\frac{i}{2Vy}F(z)-F(z)\frac{d}{dz
}\K_{1}(z)\right)-\sum_{n=1}^\infty nb_n e(nz),\\
&=\frac{1}{2\pi
i}\left(\frac{i}{2}F(z)-F(z)\frac{d}{dz}\K_{1}^+(z)\right)-\sum_{n=1}^\infty
nb_n e(nz).\tag 5.2 \endalign
$$
Let us write the Fourier expansion of $\Pi_{hol}(f(z)\K_{1}(z))$
as
$$
\Pi_{hol}(f(z)\K_{1}(z)) =\sum\limits_{m=1}^\infty d_m e(mz).
$$
If we now substitute the Fourier expansions of $F$ and $\K_{1}$ into (5.2) we find the formula
$$
d_m=-mb_m+\frac{a_m}{4\pi m}-\sum_{l=1}^{m-1}\frac{a_l}{l} (m-l)k(m-l).
$$
However, from Theorem 3.2 and Proposition 5.1,  for all $m
\geqslant 1$, we also have that
$$
d_m=mL^{++}_m(1)+\frac{a_m}{4\pi m}+\sum_{l=1}^{m-1}a_l k(m-l)+a_m\left(K +\frac{\gamma+\log 4\pi m}{V}\right).
$$
Therefore, by the definition of the Dirichlet series $L^{+}_{m}$,
as first stated in Theorem 1.2, we conclude that for all $m \geqslant
1$, we have
$$
b_m=-L^{+}_m(1). \tag 5.3
$$
Substituting (5.3) into (5.2) yields the Fourier expansion claimed
in Theorem 1.3, whose proof is now complete.

\vskip .25in \head \bf \S 6. The Dirichlet series $L_{m}^{-}$ and
the proofs of Theorems 1.4 and 1.5
\endhead

\vskip .10in Let us first prove Theorem 1.4.  To do so, we start
with (2.2) and, using Corollary 2.2, obtain the formula
$$
\frac{d}{dz}\overline{\K_2}(z)=\frac{-i}{2Vy}\overline{F}(z)+\frac{2\pi i}
V \sum_{n=1}^\infty \overline{a_n} \overline{W^*}(nz) +\frac{d}{dz}\overline{B(z)}
+\overline{F}(z)\frac{d}{dz}\overline{\K_1}(z) \tag 6.1
$$
which, in particular, implies that
$$
\Pi_{hol}\left(\frac{d}{dz}\overline{\K_2}(z)\right)=\Pi_{hol}\left(\frac{-i}{2Vy}
\overline{F}(z)+\frac{2\pi i}V \sum_{n=1}^\infty \overline{a_n} \overline{W^*}(nz)
+\frac{d}{dz}\overline{B(z)}+\overline{F}(z)\frac{d}{dz}\overline{\K_1}(z)\right).\tag 6.2
$$
By Corollary 3.4, the left-hand-side of (6.2) is zero.  Using
Theorem 3.2 and equations (3.2) and (3.3), we can compute the
$m$-th Fourier coefficient of the right-hand-side, which, since
the left-hand-side vanishes, is necessarily zero.  That is, we
have that
$$
0 = 4\pi m \int_0^\infty \int_0^1
\left(\frac{-i\overline{F}(z)}{2Vy}+\frac{2\pi i}V
\sum_{n=1}^\infty \overline{a_n} \overline{W^*}(nz)
+\frac{d}{dz}\overline{B(z)}+\overline{F}(z)\frac{d}{dz}\overline{\K_1}(z)\right)
y^{s-1}\overline{e(mz)}\, dx dy\Big|_{s=1}.
$$
In order to evaluate this,  substitute the Fourier expansions
for $F$ and $\K_{1}$, as well as the formula, 
$$
W^{\ast}(nz) = \Gamma(0,4\pi y)e^{4\pi y} e(z).
$$
Upon integrating with respect to $x$, we produce the equality
$$
\align (2\pi i)4\pi m \int_0^\infty &\left( m \overline{b_{-m}}
e^{-4\pi my} + \sum_{l=1}^\infty
\frac{\overline{a_l}}{l}(m+l)\overline{k(-m-l)}
e^{-4\pi(m+l)y}\right) y^{s-1}\, dy\\ &=\frac{2\pi i m
\G(s)}{(4\pi)^{s-1}}\left(
\frac{\overline{b_{-m}}}{m^{s-1}}+\overline{L_m^-(s)}\right),
\endalign
$$
where
$$
L_{m}^{-}(s) = \sum_{n=1}^\infty \frac{a_n}{n}\frac{
k(-m-n)}{(m+n)^{s-1}}.
$$
Now, by taking $s=1$, we get that
$$
b_{-m}=-L_m^-(1) \tag 6.3
$$
for all $-m < 0$ or $m > 0$, provided, of course, that $L_m^-(s)$
has an analytic continuation to $s=1$ which would then allow for
the above computations.  The verification that $L_{m}^{-}(s)$
admits a meromorphic continuation to, (and is analytic at), $s=1$ will be 
 completed in section 9.  In effect,
we will argue as follows. Recall that we have already used the
bounds $a_n \ll n$ and $k(n)\ll |n|^{1+\epsilon}$.  Observe that
these bounds prove $L^-_m(s)$ is
absolutely and uniformly convergent to an analytic function for
$\Re (s)> 3$.  In section 9 we will prove the functional equation
$$
L^-_m(s)=mL^-_m(s+1)+\frac{2i(4\pi)^s}{\G(s+1)}\langle
y^2f(z)\frac{d}{d\overline{z}}
\K_1(z),\overline{U_m(z,s)}\rangle,\tag 6.4
$$
where the Poincar\'e series $U_{m}(z,s)$ was introduced in section
3.  From (6.4), it is immediate that $L_{m}^{-}$ does not have a
pole at $s=1$, which then completes the proof of Theorem 1.4.

\vskip .10in Furthermore, this work yields Theorem 1.5.  
Indeed, from (5.3) and (6.3) we have shown that
$$
b_{m} =  \left\{ \matrix \displaystyle -L_{m}^{+}(1) & m \geqslant 1 \\  \\
\displaystyle -L_{-m}^{-}(1) & m \leqslant -1 \endmatrix \right\}.
$$
Substituting into (2.2), and using the equations (2.3) and (2.4),
then completes the proof of Theorem 1.5.

\vskip .10in {\bf Remark 6.1.} As an aside, let us study the
right-hand-side of (6.2) and show that it can be reduced further.
Let $H(z)=\overline{F}(z)\frac{d}{dz}\overline{\K_1}(z)$, so then,
when using the relation $F(\g z)=F(z)+\langle \g,f\rangle$, we
have that
$$
H(\g z)=j(\g,z)^2\left( H(z)+\overline{\langle \g,f\rangle}
\frac{d}{dz}\overline{\K_1}(z)\right).\tag 6.5
$$
In other words, $H(z)$ is a weight two, second-order automorphic
form. For any $g \in S_{2}(\G)$, we claim that $\langle H,g
\rangle$ is well-defined.  To see this, first choose a fundamental
domain $\F$ for $\G\backslash\H$ and, for now, let
$$
\langle H,g \rangle_{\F} =\int_\F y^2 H(z) \overline{g}(z)\, d\mu
(z).
$$
For any $\g \in \G$, it is easy to show, using the transformation
property for $g$ and (6.5), that
$$
\langle H,g \rangle_{\gamma \F} = \int_{\g \F} y^2 H(z)
\overline{g}(z)\, d\mu (z)=\int_\F y^2 H(z)\overline{g}(z)\, d\mu
(z) + \overline{\langle \g,f\rangle}\int_\F y^2
\frac{d}{dz}\overline{\K_1}(z)\overline{g}(z)\, d\mu (z).
$$
By Corollary 3.4,
$$
\langle \frac{d}{dz}\overline{\K_1},g \rangle=0,
$$
which shows that $\langle H, g \rangle_{\F}$ is $\G$ invariant,
hence $\langle H, g \rangle $ is well-defined as claimed.
Consequently, $\Pi_{hol}(H)$ makes sense and hence exists. Similar
reasoning applies to the remaining part on the right-hand-side of
(6.2).  As a result, since
$$
\Pi_{hol}\left(\frac{d}{dz}\overline{\K_2}(z)\right) = 0,
$$
by Corollary 3.4, (6.2) can be written as
$$
0=\Pi_{hol}\left(-\frac{i}{2Vy}\overline{F}(z)+\frac{2\pi i}V
\sum_{n=1}^\infty \overline{a_n} \overline{W^*}(nz)
+\frac{d}{dz}\overline{B^-(z)}\right)
+\Pi_{hol}\left(\overline{F}(z)\frac{d}{dz}\overline{\K_1}(z)\right).
$$
Possible implications of this identity have not been investigated
here.

\vskip .25in \head \bf\S 7. Bounding the Fourier coefficients of
$\K_{1}$ and $\K_{2}$
\endhead

\vskip .10in  In this section we estimate the size of the Fourier
coefficients of $E(z,s)$ and $E^*(z,s)$.  The calculations are
used to bound $k(n)$, $L_m^+(1)$, and $L_m^-(1)$.  

\vskip .10in To begin, we need the following general result.

\vskip .10in \proclaim{Lemma 7.1} Suppose $D(z)=D(x+iy)$ is a smooth
function on $\H$ which is $\Gamma$ invariant.  Assume there is a
continuous function $B(y)$ such that for each cusp $\ca$ of
$\Gamma$, we have that $|D(\sa z)| \leqslant B(y)$ as $y \rightarrow
\infty$. Then we also have
$$
\align
D(\sa z)&\ll 1 \,\,\,\,\, \text{\rm as $y\rightarrow
0$, if $B$ is decreasing and}\\
D(\sa z)&\ll B(C/y) \,\,\,\,\, \text{\rm as $y\rightarrow
0$, if $B$ is increasing }
\endalign
$$
where both implied constants and $C>0$ depend only on $D$ and $\G$ (and are independent of $x$).

\p By conjugation we may assume (as we have been doing all along) that $\ci$ is a cusp of $\GH$ and
that $\G_\ci$ is generated by the translation $z \mapsto z+1$. Let
$\Cal F_\ci=\big\{z\in \H \; \big | \; |\Re (z)|\leqslant \frac12
\big\}$, and let $\Cal F$ be the (Ford) fundamental domain for
$\GH$ defined by  $\Cal F=\big\{z\in \Cal F_\ci \; \big | \;
1<|j(\g,z)| \text{ for all } \g \in \G-\G_\ci\big\}$.  

The first statement of the Lemma is easy to prove. If $B$ is decreasing then, since $D$ is smooth, $D(z)$ is bounded on $\Cal F$ and hence on $\H$ since it is $\G$ invariant. The bound $D(\sa z)\ll 1 $ as $y\rightarrow 0$ follows trivially.

Let us assume now that $B(y)$ is increasing as $y \rightarrow \infty$. For the cusp at infinity, by assumption, 
$$
D(w) \ll B(\Im (w)) \,\,\,\,\, \text{\rm as }\Im (w)\rightarrow \infty.
$$
We next consider what happens as $w \in \Cal F$ approaches a cusp $\ca \in \R$. Set $w= \sa w'$ so that $w \rightarrow \ca$ as $\Im(w') \rightarrow \infty$. It is easy to check that
$$
\frac 1{\Im (w')} \ll \Im (\sa w') \ll \frac 1{\Im (w')} \,\,\,\,\, \text{\rm as }\Im (w')\rightarrow \infty
$$
if $\sa$ is not upper triangular and that
$$
\Im (w') \ll \Im (\sa w') \ll \Im (w') \,\,\,\,\, \text{\rm as }\Im (w')\rightarrow \infty
$$
if $\sa$ is upper triangular. Since $\Im (w) \rightarrow 0$ as $\Im (w') \rightarrow \infty$ it must be the case that $\sa$ is not upper triangular and hence, for some $C>0$,
$$
\Im (w') \leqslant \frac C{\Im (w)}. 
$$
By assumption
$$
D(\sa w') \ll B( \Im (w')) \,\,\,\,\, \text{\rm as }\Im (w')\rightarrow \infty.
$$
Therefore 
$$
D(w) \ll B( C/\Im (w)) \,\,\,\,\, \text{\rm as }w\rightarrow \ca
$$
in $\Cal F$ and it follows that, for a possibly larger $C$,
$$
D(w) \ll B(\Im (w)+ C/\Im (w))
$$
for all $w \in \Cal F$.

Now, for any  $z \in \Cal F_\ci -\Cal F$, there exists $\g \in
\G-\G_\ci$ such that $\g z=w \in \Cal F$.  It can be show that
$$
y \leqslant \Im (w) \ll \frac{1}{y},
$$
where the first inequality comes from the definition of $\Cal F$
and the second from Lemma 1.25 of \cite{Sh} (see also Proposition
2.5 of \cite{G-O'S}).  The implied constant in the upper bound
depends only on $\G$.  It now follows that, for any $z$ in $\H$, we have
$$
D(z) \ll B(C/y) \text{ \ \ as \ \ }y \rightarrow 0
$$
with the implied constant and (larger) $C$ depending only on $D$ and $\G$.   Finally, to prove the same bound for $D(\sa
z)$ we may use the same proof applied to the conjugate group $\G'=\sa^{-1}\G \sa$. Specifically, let $D'(z)=D(\sa z)$ then $D'$ is a smooth $\G'$ invariant function. Now if $\ca, \cb, \cc, \dots$ are a set of inequivalent cusps for $\GH$ then $\ca'=\sa^{-1}\ca=\ci, \cb'=\sa^{-1}\cb, \cc'=\sa^{-1}\cc, \dots$ are a set of inequivalent cusps for $\G' \backslash \H$ with corresponding scaling matrices
$\sigma_{\ca'}=\sa^{-1}\sa, \sigma_{\cb'}=\sa^{-1}\sb, \sigma_{\cc'}=\sa^{-1}\sc, \dots$. Therefore, for any cusp $\cb'$ of $\G'$ we have
$$
|D'(\sigma_{\cb'} z)|=|D(\sb z)| \leqslant B(y) \text{ \ \ as \ \ }y \rightarrow \infty.
$$
It now follows from our previous work that $ D'(z) \ll B(C/y) $ as $y \rightarrow 0$, completing the proof.
\bs

\vskip .10in To continue, we recall that equation (6.19) of
\cite{Iw1} states an explicit bound for the Fourier coefficients
of the first-order Eisenstein series, namely
$$
\phi_{\ca\cb}(n,s) \ll |n|^\sigma+|n|^{1-\sigma},
$$
with an implied constant depending on $s$ and $\G$.  We will prove
our stated bounds for the Fourier coefficients of $\K_{2}$ by
making this bound for $\phi_{\ca\cb}(n,s)$ more precise, as well as extend
the result to the functions $\phi^{*}_{\ca\cb}(n,s)$.  The main
technical result of this section is the following.

\vskip .10in \proclaim{Proposition 7.2} For each compact set $S$
in $\C$ there exist smooth functions $\psi_1(s), \psi_2(s)$ and
holomorphic functions $\xi_1(s), \xi_2(s)$ so that
$$
\align
|\phi_{\ca \cb}(n,s)| & \leqslant \frac{\psi_1(s)}{|\xi_1(s)|}
\left( |n|^\sigma+|n|^{1-\sigma} \right),\tag 7.1\\
|\phi^*_{\ca \cb}(n,s)| & \leqslant \frac{\psi_2(s)}{|\xi_2(s)|}
(\log |n|+1)\left( |n|^\sigma+|n|^{1-\sigma} \right)\tag 7.2
\endalign
$$
for all $s$ in $S$ and all $n \neq 0$. The functions $\psi_1,
\xi_1$ depend on $S$ and $\G$, and the functions 
$\psi_2, \xi_2$ depend on $S$, $\G$ and $f$.

\vskip .10in \p The bound (7.1) will follow from the proof of the
meromorphic continuation of the first-order Eisenstein series
$E_{\ca}(z,s)$ as given in Proposition 6.1 of \cite{Iw1}.  After
proving (7.1), we then employ the same method of proof, this time
using the meromorphic continuation of the second-order Eisenstein
series $E^{\ast}_{\ca}(z,s)$ as given in Theorem 3.8 of
\cite{O'S1}.  For ease of notation, $\psi$ and $\xi$ will always
represent smooth and holomorphic functions respectively, though
the functions themselves may change from line to line.

\vskip .10in From Proposition 6.1 of \cite{Iw1}, we have the
following (weaker) form of the stated result.  Given a compact
subset $S$ of $\C$, there exist functions $A_\ca(s) \not\equiv 0$
on $S$ and $A_\ca(z,s)$ on $\H \times S$ such that:

\vskip .05in
\roster
\item $A_\ca (z,s)=A_\ca (s) E_\ca (z,s)$ on $\{s \ | \ \Re (s) >1\} \cap S$,
\item $A_\ca(s)$ and $A_\ca (z,s)$ are holomorphic in $s$,
\item $A_\ca (\sb z,s) \ll e^{\varepsilon y}$ for each cusp $\cb$ and $y\geqslant 1$ say.
The implied constant depends on $\varepsilon>0, s$ and $\G$.
\endroster

\vskip .05in \noindent Furthermore, from the proof of Proposition
6.1 in \cite{Iw1}, specifically (6.1), we conclude there exists a
smooth function $\psi(s)$ on $S$ so that
$$
A_\ca (\sb z,s) \operatornamewithlimits{\ll}_{\varepsilon, \G}
\psi(s)\ e^{\varepsilon y}.
$$
The Fourier expansion of $E_{\ca}(z,s)$, namely
$$
E_\ca(\sb z,s)=\delta_{\ca \cb}y^s +\phi_{\ca
\cb}(s)y^{1-s}+\sum_{n\neq 0} \phi_{\ca \cb}(n,s)W_s(nz)
$$
then gives
$$
\phi_{\ca \cb}(n,s)\cdot 2\sqrt{|n|y} K_{s-1/2}(2\pi |n|y)
=\int_0^1 E_\ca(\sb z,s)e^{-2\pi inx} \,dx.
$$
Using (1) above and Lemma 7.1, we then obtain the bound
$$
\phi_{\ca \cb}(n,s)\cdot 2\sqrt{|n|y} K_{s-1/2}(2\pi
|n|y) \ll \frac{\psi(s)}{|A_\ca(s)|}e^{\varepsilon/y}\tag
7.3
$$
as $y\rightarrow 0$.  The $K$-Bessel function can be bounded using
the estimate
$$
\sqrt{y} K_{s-1/2}(y)=\sqrt{\frac \pi 2}e^{-y}
\left(1+O\left(\frac{1+|s|^2}{y}\right)\right) \tag 7.4
$$
for $y>1+|s|^2$, which we quote from \cite{Iw1}, formula B.36.
With this, and upon setting $y=1/\sqrt{|n|}$, we get the auxiliary
estimate
$$
 \phi_{\ca \cb}(n,s) \ll
\frac{\psi(s)}{|A_\ca(s)|}e^{3\pi\sqrt{|n|}}.
$$
Consequently,
$$
\align \left| E_\ca(\sb z,s)\right| &=\left|\delta_{\ca \cb}y^s
+\phi_{\ca \cb}(s)y^{1-s}+\sum_{n\neq 0}
\phi_{\ca \cb}(n,s)W_s(nz)\right|\\
&\ll y^\sigma +|\phi_{\ca \cb}(s)|y^{1-\sigma}+\frac{\psi(s)}{|\xi(s)|}
\sum_{n\neq 0}e^{3\pi \sqrt{|n|}-2\pi|n|y}\\
&\ll \frac{\psi(s)}{|\xi(s)|}\left(y^\sigma +y^{1-\sigma}+e^{2\pi (1/2-y)}\right)
\endalign
$$
for $y\geqslant 1$ say.  Now, repeat the argument yielding (7.3) with
this new bound to get
$$
 \phi_{\ca \cb}(n,s) \ll \frac{\psi(s)}{|\xi(s)|}
\frac{y^{-\sigma}+y^{-(1-\sigma)}} {\sqrt{|n|y}K_{s-1/2}(2\pi |n|
y)}.\tag 7.5
$$
Letting $y=1/|n|$, the proof of (7.1) is complete.

\vskip .10in An easy consequence of (7.1) that we shall need shortly is the next result.

\vskip .10in\proclaim{Corollary 7.3} For each compact set $S$ in
$\C$ there exist $\psi$ smooth and $\xi$ holomorphic such that
$$
E_\ca(\sb z,s)-\delta_{\ca \cb}y^s-\phi_{\ca \cb}(s)y^{1-s} \ll
\frac{\psi(s)}{|\xi(s)|}e^{-2\pi y}
$$
as $y\rightarrow \infty$ and
$$
E_\ca(\sb z,s)\ll \frac{\psi(s)}{|\xi(s)|} \left( y^\sigma +
y^{-\sigma} + y^{1-\sigma} + y^{\sigma-1}\right)\tag 7.6
$$
for all $y$ in $(0,\infty)$ and all $s\in S$. The implied
constants depending only on $S$ and $\G$.
\endproclaim

\vskip .10in The proof of (7.2) follows the same pattern, in this
case using Theorem 3.8 of \cite{O'S1} rather than Proposition 6.1
of \cite{Iw1}.  To begin, for any compact $S \subset \C$, there
are functions $A^*_\ca(s) \not\equiv 0$ on $S$ and $A^*_\ca(z,s)$
on $\H \times S$ such that:

\vskip .05in
\roster
\item $A^*_\ca (z,s)=A^*_\ca (s) E^*_\ca (z,s)$ on $\{s \ | \ \Re (s) >2\} \cap S$,
\item $A^*_\ca(s)$ and $A^*_\ca (z,s)$ are holomorphic in $s$,
\item $A^*_\ca (\sb z,s) \ll e^{\varepsilon y}$ for each cusp
$\cb$, $y\geqslant 1$ and implied constant depending on $\varepsilon>0, s,f$ and $\G$.
\endroster

\vskip .05in \noindent Following the method of proof of
Proposition 6.1 in \cite{Iw1}, the analysis in \cite{O'S1} yields
the bound
$$
A^*_\ca (\sb z,s) \operatornamewithlimits{\ll}_{\varepsilon,f,
\G}\psi(s)e^{\varepsilon y}.
$$
Lemma 7.1 applies to a weight zero function ($\Gamma$ invariant).
For this, we study
$$
G_\ca(z,s)=E^*_\ca(z,s)+F_\ca(z) E_\ca(z,s),
$$
which, as stated in section 1, is $\Gamma$ invariant.  Let us
write
$$
F_\ca(\sb z)=2\pi i \int_\ca^{\sb z} f(w)\, dw=2\pi i
\int_{\sb^{-1}\ca}^{ z} g(w)\, dw
$$
with $g(z)=f(\sb z)/j(\sb,z)^2 \in S_2(\sb^{-1}\G \sb)$.
Therefore,
$$
\int_z^{z+1} g(w)\,dw=0 \,\,\,\,\,\text{\rm and} \,\,\,\,\,
g(z)\ll 1/y
$$
(see (5.3), \cite{Iw2}). Consequently, we have, for each pair of cusps $\ca$,
$\cb$ and all $y$ in $(0,\infty)$, the bound
$$
F_\ca(\sb z)\operatornamewithlimits{\ll}_{f, \G} |\log y|+1. \tag 7.7
$$
(Note: This estimate improves Lemma 1.1 of \cite{O'S1}; see also
\cite{Ri}, \cite{P-R} for a different approach to this and similar bounds.) Continuing, the
bounds for the Eisenstein series $E^{\ast}(z,s)$ and $E(z,s)$,
together with (7.7) imply that as $y
\rightarrow \infty$, we have
$$
G_\ca(\sb z,s)  \ll
\frac{\psi(s)}{|\xi(s)|}e^{\varepsilon y}.
$$
Thus by Lemma 7.1,
$$
 G_\ca(\sb z,s)  \ll
\frac{\psi(s)}{|\xi(s)|}e^{\varepsilon /y} \,\,\,\,\,\text{\rm as
$y \rightarrow 0$.}
$$
With (7.6) and (7.7), we then obtain
$$
E^*_\ca(\sb z,s) \ll
\frac{\psi(s)}{|\xi(s)|}e^{\varepsilon /y} \,\,\,\,\,\text{\rm as
$y \rightarrow 0$.}
$$
By repeating the argument used to prove (7.1), we get the
auxiliary estimate
$$
 \phi^*_{\ca \cb}(n,s) \ll
\frac{\psi(s)}{|\xi(s)|}e^{3\pi\sqrt{|n|}}
$$
so then
$$
E^*_\ca(\sb z,s) \ll \frac{\psi(s)}{|\xi(s)|}y^{1-\sigma}
\,\,\,\,\,\text{\rm as $y \rightarrow \infty$.}
$$
Therefore
$$
\align
G_\ca(\sb z,s) &\ll \frac{\psi(s)}{|\xi(s)|}(|\log y|+1)(y^\sigma+y^{1-\sigma})
\text{ \ \ as \ \ }y\rightarrow \infty,\\
\text{implies \ \ } G_\ca(\sb z,s) &\ll \frac{\psi(s)}{|\xi(s)|}(|\log y|+1)
(y^{-\sigma}+y^{\sigma-1}) \text{ \ \ as \ \ }y\rightarrow 0,\\
\text{implies \ \ } E^*_\ca(\sb z,s) &\ll \frac{\psi(s)}{|\xi(s)|}(|\log y|+1)
(y^{-\sigma}+y^{\sigma-1}) \text{ \ \ as \ \ }y\rightarrow 0.
\endalign
$$
With this improved bound the equality $(7.2)$ follows in the same
manner that $(7.1)$ was proved. This completes the proof of
Proposition 7.2. \bs

\vskip .10in The analogue of Corollary 7.3 follows from Proposition 7.2.

\proclaim{Corollary 7.4} For $s$ contained in a compact set $S$ in
$\C$ we have $\psi$ smooth and $\xi$ holomorphic with
$$
E^*_\ca(\sb z,s)-\phi^*_{\ca \cb}(0,s)y^{1-s} \ll \frac{\psi(s)}{|\xi(s)|}e^{-2\pi y}
$$
as $y\rightarrow \infty$ and the implied constant depending only on $S,f$ and $\G$.
\endproclaim

\vskip .10in Another
consequence of Proposition 7.2 gives our desired bounds for the sequence $\{ b_n \}$.

\proclaim{Corollary 7.5} For every $n$, write
$$
\phi^*_{\ca \cb}(n,s)=\sum_{m=-1}^\infty b_{\ca \cb}(n,m) (s-1)^m.
$$
Then for every $m \geqslant -1$ and every $\epsilon>0$ we have
$$
b_{\ca \cb}(n,m) \operatornamewithlimits{\ll}_{m,\epsilon,f,\G} |n|^{1+\epsilon}.
$$
In particular $b_n = b(n,0) \ll |n|^{1+\epsilon}$. \p Let
$C_\epsilon$ be a circular loop around 1 with small radius
$\epsilon$. We know that
$$
b_{\ca \cb}(n,m)=\frac{1}{2\pi i} \int_{C_\epsilon}
\frac{\phi^*_{\ca \cb}(n,s)}{(s-1)^{m+1}}\, ds
$$
and by (7.2) the desired conclusion follows. \bs

\vskip .10in \flushpar {\bf Bounding the coefficients $k(n)$ of $\K_1$: }
The proof of Corollary 7.5 also applies directly to the definition of
$\K_{1}$ (on replacing (7.2) with (7.1)) to give the bounds
$$
k(n),k(-n) \ll n^{1+\epsilon}
$$
for any $\epsilon > 0$, as asserted in Theorem 1.1. The identity
$k(n)=\overline{k(-n)}$ follows from the symmetry
$\overline{E(z,\overline{s})}=E(z,s)$ because
$\overline{W_{\overline{s}}(z)}=W_s(-\overline{z})$ and therefore
$\phi_{-n}(s)=\overline{\phi_n}(\overline{s})$. With this, the
proof of Theorem 1.1 is complete.\bs

\vskip .10in \flushpar {\bf Bounding the Fourier coefficients of
$\K_{2}$ :}  The Fourier coefficients of $\K_{2}$ are expressed in
terms of the Fourier coefficients of $F$, $\K_{1}$, and the
sequence $\{b_{n}\}$.  Known results bound the Fourier
coefficients of $F$, Theorem 1.1 (whose proof is now complete)
bounds the Fourier coefficients of $\K_{1}$, and Corollary 7.5
bounds the elements of the sequence $\{b_{n}\}$.  Though it
remains to prove that $b_{m} = L_{m}^{+}(1)$ for $m \geqslant 1$ and
$b_{m} = L_{-m}^{-}(1)$ for $m \leqslant -1$, the Fourier coefficient
bounds are complete nonetheless.  To continue, let us further
analyze the Fourier coefficients $\{b_{m}\}$.

\vskip .10in If $f$ has Fourier coefficients $\{a_n\}$ in $\R$ for
all $n>0$ then we want to show that the Fourier coefficients of
$\K_{2}$ $\{b_m\}$ are also in $\R$,
provided we have $\iota(\G)=\G$ for $ \left(\smallmatrix a & b \\
c & d
\endsmallmatrix\right) \overset \iota \to \longrightarrow
\left(\smallmatrix -a & b \\ c & -d \endsmallmatrix\right). $ The
map $\iota$  is an automorphism  of $\PSL_2(\R)$, and it is easily
verified that $\g(-\overline z)=-\overline{(\iota (\g) z)}$ for
any $\g \in \PSL_2(\R)$. From this it follows that $E(-\overline
z,s)=E(z,s)$ for any subgroup $\G$ of $\PSL_2(\R)$ with
$\iota(\G)=\G$, and hence $\phi_m(s)=\phi_{-m}(s)$. Since $f$ has
real Fourier coefficients we see that
$$
\langle \iota(\g),f \rangle = \overline{\langle \g,f \rangle},
$$
and then
$$
\align E^*(z,s)&=\sum_{\g \in \G_\infty\backslash\G} \langle \g,f
\rangle\Im(\g z)^s
=\sum_{\g \in \G_\infty\backslash\G} \langle \iota(\g),f \rangle\Im(\iota(\g) z)^s\\
&=\sum_{\g \in \G_\infty\backslash\G} \overline{\langle \g,f
\rangle}\Im(\g (-\overline{z}))^s
=\overline{E^*(-\overline{z},\overline{s})}.
\endalign
$$
Therefore, $\phi^*_{m}(s)=\overline{\phi^*_m}(\overline{s})$ which
implies that $b_m=\overline{b_m}$.

\vskip .25in \head \bf \S 8. Poincar\'e series: Proofs of Theorems
3.1 and 3.2 \endhead

\vskip .10in We now prove Theorem 3.1 and Theorem 3.2.  In
essence, the material in this section is based on \cite{Se},
Chapter 17 of \cite{Iw3} and \cite{Ne}.  The weight $k$ Poincar\'e
series is defined by the series
$$
P_{\ca m}(z)_k=\sum_{\g \in \G_\ca\backslash\G}
\frac{e(m\sa^{-1}\g z)}{j(\sa^{-1}\g, z)^k}. \tag 8.1
$$
The series (8.1) converges absolutely and uniformly if $k > 2$ but
not when $k=2$.  Hecke addressed this problem by introducing a
complex parameter $s$ and taking a limit.  We will follow this
approach employing the non-holomorphic Poincar\'e series
$U_{\ca m}(z,s)$ from section 3. If $m=0$ we have that $U_{\ca
0}(z,s)=E_{\ca}(z,s)$.  Since the non-holomorphic Eisenstein
series is absolutely convergent for $\Re (s)>1$, we have that the function $E_{\ca}(z,\Re (s))$ is a majorant of
$U_{\ca m}(z,s)$ for
$m\geqslant 0$.

\vskip .10in \proclaim{Lemma 8.1} For $m\geqslant 1$ and $\Re
(s)>1$ the Poincar\'e series $U_{\ca m}(z,s)$ is square
integrable, i.e. $U_{\ca m}(z,s)$ is in $L^2(\G \backslash \H)$.
\p We first examine the size of $U_{\ca m}$ in the neighborhood of
each cusp. Setting $s=\sigma +it$, we have
$$
\align \left|U_{\ca m}(\sa z,s)\right| &\ll y^\sigma e^{-2\pi
my}+\sum \Sb \g \in \G_\ca\backslash\G
\\ \g \neq identity \endSb \Im(\sa^{-1}\g\sa z)^\sigma \\
&\ll y^\sigma e^{-2\pi my}+|E_\ca(\sa z,\sigma)-y^\sigma| \ll 1.
\endalign
$$
At any other cusp $\cb \neq \ca$
$$
U_{\ca m}(\sb z,s)\ll E_\ca(\sb z,\sigma)\ll \phi_{\ca \cb}(s)y^{1-\sigma} \ll 1.
$$
In other words, $U_{\ca m}$ is bounded on $\G \backslash \H$ and
hence in $L^2(\G \backslash \H)$ since $\G \backslash \H$ has
finite volume. \bs

\vskip .10in We will study the Poincar\'e series $U_{\ca m}(\sa
z,s)$ by means of its spectral expansion, which we now recall
(see, for example, \cite{Iw1} and references therein for further
background information and complete proofs). The hyperbolic
Laplacian $\Delta$ operates on the space $L^2(\G \backslash \H)$,
and any element $\xi$ of $L^2(\G \backslash \H)$ may be decomposed
into constituent parts from the discrete and continuous spectrum
of $\Delta$. This decomposition, often referred to as the
Roelcke-Selberg expansion, amounts to the identity
$$
\xi(z)=\sum_{j=0}^\infty\langle \xi,\eta_j\rangle
\eta_j(z)+\frac{1}{4\pi }\sum_\cb \int_{-\infty}^{\infty}\langle
\xi ,E_\cb(\cdot,1/2+ir)\rangle E_\cb(z,1/2+ir)\,dr,\tag 8.2
$$
where $\{\eta_j\}$ denotes a complete orthonormal basis of Maass
forms, with corresponding eigenvalues $\lambda_j=s_j(1-s_j)$,
which forms the discrete spectrum.  For notational convenience, we
wrote $\langle \cdot, \cdot \rangle=\langle \cdot, \cdot
\rangle_0$ for the inner product on $\Gamma \backslash \H$ of
weight zero forms (i.e. $\G$-invariant functions).  As always, we
will write $s_j=\sigma_j+it_j$, chosen so that $\sigma_j \geqslant
1/2$ and $t_j \geqslant 0$, and we enumerate the eigenvalues,
counted with multiplicity, by $0=\lambda_0< \lambda_1 \leqslant
\lambda_2 \leqslant \dots $. For each $j$, the Fourier expansion
of $\eta_j$ is
$$
\eta_j(\sa z)=\rho_{\ca j}(0)y^{1-s_j}+\sum_{m\neq 0}
\rho_{\ca j}(m)W_{s_j}(mz).\tag 8.3
$$
For all but finitely many of the $j$ (corresponding to
$\lambda_{j} < 1/4$) we have $\sigma_j=1/2$ and $\rho_{\ca
j}(0)=0$.   The expansion (8.2) is absolutely convergent for each fixed $z$ and
uniform on compact subsets of $\H$, provided $\xi$ and $\Delta
\xi$ are smooth and bounded (see, for example, Theorem 4.7 and
Theorem 7.3 of \cite{Iw1}).  By taking $\xi = U_{\ca m}$, we then
obtain the spectral expansion for the Poincar\'e series, which
yields the identity
$$
\multline
U_{\ca m}(z,s)\pi^{-1/2}(4\pi m)^{s-1/2}\G(s)=\sum_{j=1}^\infty \G(s-s_j)
\G(s-1+s_j)\overline{\rho_{\ca j}}(m)\eta_j(z)\\
+\frac{1}{4\pi}\sum_\cb \int_{-\infty}^\infty \G(s-1/2-ir)\G(s-1/2+ir)
\overline{\phi_{\ca \cb}}(m,1/2+ir)E_\cb(z,1/2+ir)\, dr.
\endmultline \tag 8.4
$$
The expansion (8.4) includes the identity
$$
\langle U_{\ca m}(\cdot ,s),\eta_j\rangle =
\frac{\pi^{1/2}\G(s-s_j)
\G(s-1+s_j)}{(4sm)^{s-1/2}\G(s)}\overline{\rho_{\ca j}}(m),
$$
with a similar formula which evaluates the inner product of the
Poincar\'e series $U_{\ca m}(z,s)$ with the Eisenstein series
$E_{\ca}(z,s)$.  The proofs of these formulas come from unfolding
the integrals under study and unfolding the series which defines
the Poincar\'e series.  These calculations we leave for the
interested reader.  When looking toward Theorem 3.2, the
appearance of the coefficients $\rho_{\ca j}(m)$ and $\phi_{\ca
\cb}(m,1/2+ir)$ is natural since $U_{\ca m}$ isolates
$m$-th Fourier coefficients (see Theorem 3.2, and, more
specifically, see \cite{Ne} or Chapter 17 of \cite{Iw3}).

\vskip .10in Initially, (8.3) is valid for $\Re (s)>1$.  The
remainder of this section shows that (8.3) converges absolutely
and uniformly in $s$ in compact subsets not containing a number of
the form $s_j-n$ or $1-s_j-n$ for $n \in \N$. These points are
poles caused by the factors $\G(s-s_j)\G(s-1+s_j)$.  Going
further, we will prove bounds regarding the growth in $z$ of
$U_{\ca m}(z,s)$ and $\frac{d}{dz}U_{\ca m}(z,s)$.  These
computations will yield the proofs of Theorem 3.1 and Theorem 3.2.

\vskip .10in To control the size of $\rho_{\ca j}(m)$ and
$\phi_{\ca \cb}(m,1/2+ir)$ we appeal to the following formula of
Bruggeman and Kuznetsov, as stated in (9.13) of \cite{Iw1}.  With
notation as above, let 
$$
N_{\ca}(T) = \sum_{|t_j|<T} |\G(s_j) \G(1-s_j) \rho_{\ca j}(m)|^2
+\frac{1}{4\pi}\sum_\cb \int_{-T}^T |\G(1/2+ir)\G(1/2-ir)\phi_{\ca
\cb}(m,1/2+ir)|^2 \, dr.
$$
Then
$$
N_{\ca}(T) =\frac{T^2}{2\pi |m|}+O(T) \,\,\,\,\,\text{\rm as $T
\rightarrow \infty$}, \tag 8.5
$$
with an implied constant which depends solely on the discrete
group $\G$. Recall that Stirling's formula states that the
classical gamma function satisfies the bound
$$
|\G(\sigma+it)| \sim \sqrt{2\pi} |t|^{\sigma-1/2} e^{-\pi|t|/2}
\,\,\,\,\,\text{\rm as $|t|\rightarrow \infty$}.
$$
For simplicity, we may assume that $T > 0$.  By combining
Stirling's formula with (8.5), we get the bounds
$$
|\rho_{\ca j}(m)|^2 \ll \frac{|t_j|^2}{|m|}e^{\pi |t_j|}, \tag 8.6
$$
and
$$
\int_{T}^{T+1} |\phi_{\ca \cb}(m,1/2+ir)|^2\, dr \ll
\frac{T^2}{|m|}e^{\pi T}. \tag 8.7
$$
We next need bounds concerning $K$-Bessel functions and Whittacker
functions.

\vskip .10in \proclaim{Lemma 8.2} For any integer $k\geqslant 0$
we have, for $\sigma>1/2-k$, the bounds
$$
\left| K_{s-1/2}(y) \right| \ll
\frac{|s|^{2k}+1}{y^{2k-1/2+\sigma}}|\G(s)| \tag 8.8
$$
and
$$
\left|\frac{d}{dy} K_{s-1/2}(y) \right| \ll
\frac{|s|^{2k+1}+1}{y^{2k+1/2+\sigma}}|\G(s)|, \tag 8.9
$$
where the implied constant depends solely on $\sigma$ and $k$ \p
First consider the case $k=0$.  From page 205 of \cite{Iw1}, we have
the expression
$$
K_{s-1/2}(y)=\frac{1}{\sqrt{\pi}} \G(s) \left( \frac y2\right)^{1/2-s}
\int_0^\infty (u^2+1)^{-s} \cos(uy)\, du\tag 8.10
$$
which is absolutely convergent for $\sigma>1/2$.  Trivially, this
gives (8.8) with $k=0$.  Next, we recall the recursive formula
$$
K_{s-1/2}(y) =\frac{2s+1}y K_{s+1/2}(y) - K_{s+3/2}(y)
$$
which comes from integrating (8.10) through integration by parts.
The recursive relation provides the inductive step by which (8.8)
follows from (8.10) for all $k \geqslant 0$.  Similarly, (8.9) follows
from (8.10) with $k=0$, and the general case is then derived using
the indentity
$$
\frac{d}{dy} K_{s-1/2}(y)=\frac{s-1/2}y K_{s-1/2}(y)-K_{s+1/2}(y).
$$
\bs

\vskip .10in Recall that
$$
W_{s}(z) = 2 y^{1/2}K_{s-1/2}(2\pi  y) e( x).
$$
Therefore, from Lemma 8.2, we see that for any $k \geqslant 0$
and $\sigma>1/2-k$, we have the bounds
$$
W_s(nz)\ll \frac{|s|^{2k}+1}{(|n|y)^{2k-1+\sigma}}|\G(s)|,\tag
8.11
$$
and
$$
 \frac{d}{dz}W_s(nz)\ll
\left(1+\frac{|s|+1}y\right) \frac{|s|^{2k}+1}
{(|n|y)^{2k-1+\sigma}}|\G(s)|,\tag 8.12
$$
where the implied constants depend solely on $\sigma$ and $k$.

\vskip .10in
 Recall the definition of $y_\G(z)$ before Theorem 3.1. The estimates (8.6), (8.11), and (8.12) now can be combined
with the Fourier expansion (8.3) to show that if $\sigma_{j} =
1/2$, then
$$
\eta_j(z)\ll y_\G(z)^{1/2}+|t_j|^{7/2}y_\G(z)^{-3/2}, \tag 8.13
$$
and
$$
y\frac{d}{dz}\eta_j(z) \ll y_\G(z)^{1/2}+|t_j|^{9/2}y_\G(z)^{-3/2}
\tag 8.14
$$
(compare, for example, with (8.3') and (8.4) of \cite{Iw1}).
Our argument at this point shows that
$$
\sum_{j=N}^\infty \G(s-s_j)\G(s-1+s_j)\overline{\rho_{\ca
j}}(m)\eta_j(z)\ll |m|^{-1/2}y_\G(z)^{1/2},
$$
where the implied constant depends on $s$ and $\G$. Clearly, the
dependence of this bound on $s$ is uniform on compact sets not
containing $s_j-n$, $1-s_j-n$ for $n \in \N$.  In other words, the
term in (8.4) associated to the discrete spectrum admits a
meromorphic continuation to all $s \in \C$ and, as claimed in
Theorem 3.1, we have the desired growth in the cusps.  It remains
to consider the integral term in (8.4).  For this, we begin with
the following proposition, which can be compared to (7.10) of
\cite{Iw1}.

\vskip .10in
\proclaim{Proposition 8.3} For any cusp $\ca$ and $z
\in \G \backslash \H$, we have the bounds
$$
\int_T^{T+1} |E_\ca( z, 1/2+ir)|^2\, dr \ll y_\G(z)T^{10},\tag
8.15
$$
and
$$
\int_T^{T+1} |y\frac{d}{dz}E_\ca( z, 1/2+ir)|^2\, dr \ll
y_\G(z)T^{12},\tag 8.16
$$
where the implied constant depends solely on $\G$. \p The proof
will follow by studying the Fourier expansion (4.3).  From the
functional equation for the scattering matrix (Theorem 6.6 of
\cite{Iw1}), we obtain the estimate
$$
\phi_{\ca \cb}(1/2+ir) \ll 1.
$$
Therefore, with (8.11),
$$
E_\ca(\sb z, 1/2+ir) \ll  \sqrt{y}+ \sum_{m\neq 0} |\phi_{\ca
\cb}(m,1/2+ir)| (|r|^{2k} +1)|\G(1/2+ir)| (|m|y)^{2k-1/2}.
$$
Consequently
$$
\align \int_T^{T+1} &|E_\ca(\sb z, 1/2+ir)|^2\, dr \ll y
+T^{2k}e^{-\pi T/2}y^{1-2k}\int_T^{T+1} \sum_{m\neq 0} |\phi_{\ca
\cb}(m,1/2+ir)| |m|^{1/2-2k}\, dr\\ &+T^{4k}e^{-\pi
T}y^{1-4k}\int_T^{T+1} \sum_{m_1\neq 0} \sum_{m_2 \neq 0}|
\phi_{\ca \cb}(m_1,1/2+ir)\phi_{\ca \cb}(m_2,1/2+ir)| |m_1
m_2|^{1/2-2k}\, dr.
\endalign
$$
This bound is valid after we have justified an interchange of
integration and summation.  The following lemma allows one to
employ the Lebesgue dominated convergence theorem, for example, to
interchange integration and summation.

\vskip .10in \proclaim{Lemma 8.4} For $r \in [T,T+1]$ we have
$$
\phi_{\ca \cb}(m,1/2+ir)\ll |m|^2
$$
for an implied constant depending on $T$, $\G$ only.
\p As in
Proposition 7.2, we write $E_\ca(z,s)$ as a quotient of
holomorphic functions $A_\ca(z,s)/A_\ca(s)$, which is valid for
$s$ in $S$ where, in this instance, $S$ is the line segment
between $1/2+iT$ and $1/2+i(T+1)$.  Theorem 6.11 of \cite{Iw1}
states that $E_{\ca}(z,s)$ has no poles on $S$, in particular, so
we may assume, after multiplying the numerator and denominator of
$A_\ca(z,s)/A_\ca(s)$ by a polynomial if necessary, that
$A_\ca(s)$ has no zeros on $S$.  As in the proof of (7.5), and
noting that $|\phi_{\ca \cb}(1/2+ir)|\leqslant 1$, we arrive at
the bound
$$
\phi_{\ca \cb}(m,1/2+ir)\ll \frac{\psi(1/2+ir)}{|A_\ca(1/2+ir)|}
\frac{y^{-1/2}} {\sqrt{|m|y}K_{ir}(2\pi |m| y)}
$$
where $\psi$ is a smooth function, and the consideration is valid
for $r\in [T,T+1]$ and $y<1$, say.  If we set $y=(\log |m|)/(2\pi
|m|)$, we get
$$
\phi_{\ca \cb}(m,1/2+ir)\ll  \frac{\sqrt{|m|}}{\log |m|K_{ir}(\log
|m|)},
$$
with an implied constant depending on $T$ and $\G$.  By using the
asymptotic (7.4), the proof of Lemma 8.4 is complete. \bs

\vskip .10in Let us now continue with the proof of Proposition
8.3.  We apply (8.7) to see that
$$
\int_T^{T+1} |E_\ca(\sb z, 1/2+ir)|^2\, dr \ll
y+y^{1-2k}T^{2k+1}+y^{1-4k}T^{4k+2}
$$
for any $k\geqslant 2$, with an implied constant depending on $k$
and $\G$.  Estimate (8.15) of the proposition now follows when
taking $k=2$.  Estimate (8.16) is proved similarly using (8.12)
instead of (8.11). With this, the proof of Proposition 8.3 is
complete.\bs

\vskip .10in We now analyze the integral in (8.4).  Using (8.7),
(8.15), and the Cauchy-Schwartz inequality we find
$$
\align &\int_{T}^{T+1}
|\G(s-1/2-ir)\G(s-1/2+ir)\overline{\phi_{\ca
\cb}}(m,1/2+ir)E_\cb(z,1/2+ir)|\, dr\\& \ll
|\G(s-1/2-iT)\G(s-1/2+iT)|  \left(
\int_{T}^{T+1}|\overline{\phi_{\ca \cb}}(m,1/2+ir)|^2 \, dr
\int_T^{T+1} |E_\ca(\sb z, 1/2+ir)|^2\, dr\right)^{1/2}\\& \ll
|(T+t)(T-t)|^{\sigma-1}T e^{-\pi|T-t|/2 -\pi|T+t|/2+\pi T/2}
|m|^{-1/2}\sqrt{y_\G(z)T^{10}}.
\endalign
$$
Thus, for $s$ in a compact set $S$, the continuous spectrum
contribution to the spectral expansion of $U_{\ca m}(z,s)$ is
absolutely and uniformly convergent, and is bounded by
$|m|^{-1/2}\sqrt{y_\G(z)}$. The meromorphic continuation of
$U_{\ca m}(z,s)$ is therefore given by (8.4) to the right of the
line of integration at $\Re (s)=1/2$. Were we to consider $s \in
\C$ to the left of $\Re (s)=1/2$, then we would express $U_{\ca
m}(z,s)$  by (8.4) together with Eisenstein series that arise when
the line of integration is crossed  (see Satz 6.6 of \cite{Ne} or
\S 6 of \cite{C-O'S}). However, we are only concerned with $s$
near $1$.  It may now be seen from (8.4) that $U_{\ca m}(z,s)$ is
holomorphic in $s$ at $s=1$.

\vskip .10in We note that Selberg was the first to prove the
meromorphic continuation of $U_{\ca m}(z,s)$, see \cite{Se}. Our
proof above shows that
$$
U_{\ca m}( z,s) \ll |m|^{-1/2}\sqrt{y_\G(z)} \tag 8.17
$$
for $\Re (s)>1/2$ with an implied constant depending on $s$.

\vskip .10in Let $U'_{\ca m}(z,s)=\frac{d}{dz}U_{\ca m}(z,s)$. By
the same arguments, using (8.14) and (8.16), we see that $U'_{\ca m}(z,s)$ also has a meromorphic continuation to all $s$ in $\C$ and satisfies
$$
yU'_{\ca m}( z,s) \ll  |m|^{-1/2}\sqrt{y_\G(z)} \tag 8.18
$$
for $\Re (s)>1/2$.  It is also true that $U'_{\ca m}(z,s)$ is
 holomorphic in $s$ at $s=1$. With all this, the proof of Theorem 3.1
is complete.

\vskip .10in For the reasons given in Section 3, we define the
holomorphic, weight two, Poincar\'e series by
$$
P_{\ca m}(z)_2=2i U'_{\ca m}(z,1) +4\pi m V_{\ca m}(z,1).\tag 8.19
$$
It is elementary to show that the right hand side of (8.19) has
weight two.  Using the series definition for $V_{\ca m}$ and the
differential equation
$$
\bigl(\Delta -s(1-s)\bigr)U_{\ca m}(z,s)=4\pi msU_{\ca m}(z,s+1),
$$
it is easy to show that $\frac{d}{d\overline{z}}P_{\ca m}(z)_2=0$,
i.e. the form $P_{\ca m}(z)_2$ is holomorphic.  Therefore, we have
the Fourier expansion
$$
j(\sb,z)^{-2}P_{\ca m}(\sb z)_2=\sum_{n \in \Z} p_\cb(n) e(nz).
$$
By adapting the proof of Lemma 8.1, one shows that
$$
j(\sb,z)^{-2}V_{\ca m}(\sb z,1) \ll y^{-1} \text{ \ as \ }
y\rightarrow \infty.\tag 8.20
$$
Using (8.18), (8.19), and (8.20), we conclude that we must have
$p_\cb(n)=0$ for $n\leqslant 0$.  Consequently $P_{\ca m}(z)_2$ is
in $S_2(\G)$ as we wanted to show. This proves the first part of
Theorem 3.2.  The remaining aspect of Theorem 3.2 follows from a
direct computation using (3.2) that we leave to the reader.

\vskip .25in \head \bf \S 9. Proofs of Proposition 3.3 and the
meromorphic continuation of $L_{m}^{+}$ and $L_{m}^{-}$ \endhead

\vskip .10in In this section we tie up the remaining `loose
ends' by completing the proof of Propostion 3.3 and the
meromorphic continuation of $L_{m}^{+}$ and $L_{m}^{-}$, as
claimed in Theorem 1.2.

\vskip .10in For $\Re (s)$ sufficiently large, $f \in S_2(\G)$ and $F=2\pi i\int f$, define the automorphic series
$$
Q_{m}(z,s;f)=\sum_{\g \in \G_\ci \backslash \G} f(\g z) \Im(\g
z)^s e(m\g z),
$$
and
$$
 Q_{m}(z,s;F)=\sum_{\g \in \G_\ci \backslash \G}
F(\g z) \Im(\g z)^s e(m\g z).
$$
Proceeding formally, if we unfold the inner product of $\K_1$ and
$Q_m(\cdot, s;F)$, we get
$$
\langle \K_1, \overline{Q_m(\cdot, s;F)}\rangle = \int_0^\infty
\int_0^1 \K_1(z)F(z)e(mz)y^{s-2}\, dx dy,
$$
which in turn can be explicitly evaluated using the Fourier
expansions of $F$ and $\K_{1}$, yielding
$$
\int_0^\infty \int_0^1 \K_1(z)F(z)e(mz)y^{s-2}\, dx dy
 =
\frac{\G(s-1)}{(4\pi)^{s-1}} L^-_m(s).
$$
As we will see, we can manipulate this inner product to obtain
(6.4), which will provide a meromorphic continuation of
$L_{m}^{-}$.

\vskip .10in From the bound (7.7), we get that
$$
F(\g z) \ll 1+|\log \Im(\g z)|+|\log \Im( z)|.
$$
By mimicking the proof of Lemma 8.1, we immediately arrive at the
following estimate.

\vskip .10in \proclaim{Lemma 9.1} For any $f\in S_2(\G)$ and
integer $m>0$, the series $Q_{ m}(z,s;F)$ is absolutely convergent
for $\Re (s)>1$.  Furthermore, if $s = \sigma+it$ with $\sigma >
1$, we have
$$
Q_{ m}( z,s;F) \ll y_\G(z)^{1-\sigma}
$$
with the implied constant depending on $s$, $f$ and $\G$ alone.
\endproclaim

\vskip .10in For the remainder of this section, we let
$C^\infty(\G\backslash\H,k)$ denote the space of smooth functions
$\psi$ on $\H$ that transform as
$$
\psi(\g z)=\ee(\g,z)^k\psi(z)
$$
for $\g$ in $\G$ and $\ee(\g,z)=j(\g,z)/|j(\g,z)|$. For example,
one element of this space is given by the series
$$
U_{\ca m}(z,s,k)=\sum_{\g \in \G_\ca \backslash\G}\Im(\sa^{-1}\g
z)^s e(m\sa^{-1} \g z) \ee(\sa^{-1}\g,z)^{-k},\tag 9.1
$$
which is the weight $k$ non-holomorphic Poincar\'e series, or, in
particular, the Eisenstein series
$$
E_{k \ca} (z,s):=U_{\ca 0}(z,s,k)
$$
in the special case when $m=0$. (Warning: It should be clear from
the context whether we mean this new notion of weight or the
previous definition of weight.) Trivially, if $\psi \in
C^\infty(\G\backslash\H,k)$ then $|\psi|$ has weight zero (in
either definition), and $\langle \cdot, \cdot
\rangle = \langle \cdot, \cdot \rangle_0$ is an inner product for
$C^\infty(\G\backslash\H,k)$. We define the \it Maass raising and
lowering operators \rm by
$$
R_k=2iy\frac{d}{dz} +\frac{k}{2}, \ L_k=-2iy\frac{d}{d\bar z}
-\frac{k}{2}.
$$
It is an elementary exercise to show that
$$
R_k: C^\infty(\G\backslash\H,k)\rightarrow C^\infty(\G\backslash\H,k+2),\ L_k:
C^\infty(\G\backslash\H,k)\rightarrow C^\infty(\G\backslash\H,k-2),
$$
and, furthermore, the hyperbolic Laplacian $\Delta$ can be
realized as
$$
\Delta =-L_{2}R_0  = -R_{-2}L_0.  \tag 9.2
$$
By direct verification we have the next lemma (see also Lemma 4.1
of \cite{C-O'S}).

\vskip .10in \proclaim{Lemma 9.2} For any $\g \in \PSL_2(\R )$ and
any smooth function $F$, let
$$
\mu(s,k,F)=F(\g z) \Im(\g z)^s e(m \g z) \ee( \g, z)^{-k}.
$$
Then
$$
\align
R_k \mu(s,k,F) &= 2i \mu(s+1,k+2,\frac{d}{dz}F)+(s+k/2)\mu(s,k+2,F)-4\pi m \mu(s+1,k+2,F),\\
L_k \mu(s,k,F) &= -2i \mu(s+1,k-2,\frac{d}{d\overline z}F)+(s-k/2)\mu(s,k-2,F).
\endalign
$$
\endproclaim

\vskip .10in Lemma 9.2 applies in the special case $F \equiv 1$ to yield
the weight $k$ non-holomorphic Poincar\'e series identities
$$
R_k U_{\ca m}(z,s,k)=(s+k/2)U_{\ca m}(z,s,k+2)-4\pi m U_{\ca
m}(z,s+1,k+2)
$$
and
$$
L_k U_{\ca m}(z,s,k)=(s-k/2)U_{\ca m}(z,s,k-2).
$$
Using this last identity, together with our established notational
conventions, we see that
$$
\align
Q_m(z,s;f)&=f(z)\sum_{\g \in \G_\ci \backslash\G}j(\g,z)^2 \Im(\g z)^s e(m \g z) \\
&=y f(z)\sum_{\g \in \G_\ci \backslash\G}\Im(\g z)^{s-1} e(m \g z) \ee(\g,z)^{2}\\
&=y f(z)U_{ m}(z,s-1,-2) =y f(z)\frac{L_0 U_{ m}(z,s-1)}{s-1}.\tag
9.3
\endalign
$$
Next, combine Lemma 9.2 (this time with $F=2\pi i \int f$ as usual) and the identity (9.2) to get
$$
(\Delta -s(1-s))Q_{ m}(z,s;F)=4\pi s Q_{ m}(z,s+1;f)+4\pi msQ_{ m}(z,s+1;F).
$$
Finally, by taking the inner product with $\K_{1}$, we get
$$
\align
 4\pi s \langle \K_1, \overline{Q_m(\cdot,
s+1;f)}\rangle+4\pi m s &\langle \K_1, \overline{Q_m(\cdot,
s+1;F)}\rangle =\langle \K_1, (\Delta
-\overline{s(1-s)})\overline{Q_m(\cdot, s;F)}\rangle
\\
&=\langle (\Delta -s(1-s))\K_1, \overline{Q_m(\cdot, s;F)}\rangle
\\
&=\langle \Delta \K_1, \overline{Q_m(\cdot,
s;F)}\rangle-s(1-s)\langle \K_1, \overline{Q_m(\cdot,
s;F)}\rangle.\tag 9.4
\endalign
$$
All calculations yielding (9.4) are correct providing the inner
products make sense and we can justify moving $\Delta$ from one
side to the other.  For example, if all functions were bounded on
$\GH$, then the manipulations are correct (see Lemma 4.1 of
\cite{Iw1}).  Unfortunately, the functions in (9.4) are not
bounded, so further analysis is required.  The following
proposition proves the bounds required to validate (9.4).

\proclaim{Proposition 9.3} Suppose $\phi_1 \in
C^\infty(\G\backslash\H,k)$ and $\phi_2 \in
C^\infty(\G\backslash\H,k+2)$. Let $A, B \in \R$ with $A+B<0$. If
$$
\phi_1(z), R_k \phi_1(z) \ll y_\G(z)^A \,\,\,\,\,\text{\rm
and}\,\,\,\,\, L_{k+2} \phi_2(z), \phi_2(z) \ll y_\G(z)^B,
$$
then
$$
\langle R_k \phi_1,  \phi_2\rangle +\langle  \phi_1,
L_{k+2}\phi_2\rangle=0.
$$
\p Let $\epsilon>0$ be such that $A+B<-\epsilon$, and choose $s
\in (1,1+\epsilon)$. Since $E_\ca( z,s) \ll y_\G(z)^s$ the inner
products in the sum
$$
 \langle R_k \phi_1,  \phi_2 E(\cdot, \overline{s})\rangle +\langle
 \phi_1,  (L_{k+2}\phi_2) E(\cdot, \overline{s})\rangle
$$
are absolutely convergent, so then we may unfold the integrals to
get
$$
\align
 \langle R_k \phi_1,  \phi_2 E(\cdot, \overline{s})\rangle &+\langle
 \phi_1,  (L_{k+2}\phi_2) E(\cdot, \overline{s})\rangle=
\int_0^\infty \int_0^1 (R_k \phi_1(z))\overline{\phi_2}(z)y^{s-2}
\, dxdy \\&+ \int_0^\infty \int_0^1
\phi_1(z)\overline{(L_{k+2}\phi_2(z)})y^{s-2} \, dxdy. \tag 9.5
\endalign
$$
It is clearer to now replace $\int_0^\infty$ with $\int_{1/ D}^D$
and then later let $D \rightarrow \infty$. With the definitions of
$R_k$ and $L_{k+2}$, $(9.5)$ becomes
$$
\align \int_{1/ D}^D \int_0^1 &\left[ \bigl( (iy \frac d{dx}
+y\frac{d}{dy})\phi_1(z)\bigr)\overline{\phi_2}(z)+\phi_1(z)(iy
\frac d{dx} +y\frac{d}{dy})\overline{\phi_2}(z)\right]y^{s-2} \,
dxdy
\\&
-\int_{1/ D}^D \int_0^1 \phi_1(z)\overline{\phi_2}(z)y^{s-2} \,
dxdy.
\endalign
$$
Now use integration by parts with respect to both $x$ and $y$.
Observing that most terms cancel, we are left with
$$
\align \int_0^1
&\left[\phi_1(x+iD)\overline{\phi_2}(x+iD)D^{s-1}-\phi_1(x+i/D)
\overline{\phi_2}(x+i/D )D^{1-s}\right] \, dx \\& -s\int_{1/ D}^D
\int_0^1 \phi_1(z)\overline{\phi_2}(z)y^{s-2} \, dxdy.
\endalign
$$
By assumption, $\phi_1(z)\overline{\phi_2}(z) \ll y_\G(z)^{A+B}$,
hence we obtain the bounds
$$
\phi_1(z)\overline{\phi_2}(z)\ll y^{A+B} \text{ \ \ as \ \ }y
\rightarrow \infty
$$
and, by Lemma 7.1,
$$
 \phi_1(z)\overline{\phi_2}(z)\ll
1 \text{ \ \ as \ \ }y \rightarrow 0,
$$
and, indeed, the asymptotics are independent of $x$. These bounds are
just enough to show that the first integral above vanishes as $D
\rightarrow \infty$. Therefore
$$
\align \langle R_k \phi_1,  \phi_2 E(\cdot, \overline{s})\rangle
+\langle  \phi_1,  L_{k+2}\phi_2 E(\cdot, \overline{s})\rangle & =
-s\int_0^\infty \int_0^1 \phi_1(z)\overline{\phi_2}(z)y^{s-2} \, dxdy \\
&=-s\langle  \phi_1,  \phi_2 E_{-2}(\cdot, \overline{s})\rangle
\endalign
$$
for the weight $-2$ Eisenstein series defined by (9.1) for $m=0$,
$k=-2$. This is valid for $s$ in $(1,1+\epsilon)$. By analytic
continuation this is true for all $s$ with $1/2<\Re (s) <1+\epsilon$
say. Finally, equating residues at $s=1$ yields the theorem because
$E_{-2}(z, s)$ is holomorphic at $s=1$. \bs

\proclaim{Corollary 9.4} Assume $\phi_1(z)$ and $\phi_2(z)$ are
smooth of weight zero with $A+B<0$, and suppose
$$
\phi_1, \ R_0 \phi_1,  \ \Delta \phi_1 \ll y_\G(z)^A, $$ and
$$
\phi_2, \ R_0 \phi_2,  \ \Delta \phi_2 \ll y_\G(z)^B.
$$
Then
$$
\langle  \Delta \phi_1,  \phi_2 \rangle=\langle  \phi_1,
\Delta\phi_2 \rangle.
$$
\p One applies Proposition 9.3 twice and uses the identity which
expresses the Laplacian in terms of the raising and lowering
operators. \bs

\vskip .10in \flushpar {\bf Proof of Proposition 3.3: }  The proof
is an immediate consequence of Proposition 9.3 when taking $k=0$
together with the definitions of the functions under
consideration. More specifically, given the functions in
Proposition 3.3, one applies Proposition 9.3 with $\phi_{1} (z) =
\varphi_{1}(z)$ and $\phi_{2}(z) = \Im (z) \cdot \varphi_{2}(z)$,
after which one then easily computes the derivatives in question. \bs

\vskip .10in \flushpar {\bf Meromorphic continuation of
$L_{m}^{-}$: }  Corollary 9.4 implies that (9.4) holds for $\Re
(s)$ sufficiently large.  Using (2.5), we then have that
$$
\langle \Delta \K_1, \overline{Q_m(\cdot, s;F)}\rangle=\langle
-V^{-1},\overline{Q_m(\cdot, s;F)}\rangle=0 \tag 9.6
$$
where the last equality comes from unfolding the integral in
question and using that $f$ is a holomorphic cusp form.  If we now
combine (9.3), (9.4) and (9.6), we then get
$$
\frac{\G(s+1)}{(4\pi)^{s-1}}L^-_m(s)=m
\frac{\G(s+1)}{(4\pi)^{s-1}}L^-_m(s+1)+4\pi\langle \K_1,
\overline{y f(z)L_0 U_m(z,s)}\rangle.\tag9.7
$$
However, the structure of the operators $L$ and $R$ are such that
we have the relation
$$
\langle  \K_1, \overline{y f(z)L_0 U_m(z,s)}\rangle=\langle y
f(z)\K_1, R_0\overline{U_m(z,s)}\rangle=-\langle  y f(z)L_0\K_1,
\overline{U_m(z,s)}\rangle.
$$
Substituting this into (9.7) completes the proof of the identity
$$
L^-_m(s)=mL^-_m(s+1)+\frac{2i(4\pi)^s}{\G(s+1)}\langle
y^2f(z)\frac{d}{d\overline{z}}
\K_1(z),\overline{U_m(z,s)}\rangle.\tag9.8
$$
We see that $Q_m(z,s;F)$ does not appear in (9.8) and a second proof of (9.8) is to simply unfold the inner product on the right side.
The bounds on the Fourier coefficients $\{a_{n}\}$ and $\{k(n)\}$
are such that the Dirichlet series which defines $L_{m}^{-}(s)$
converges for $\Re (s) > 3$.  The bound (8.17) and identity (9.8) provide
the meromorphic continuation to $\Re (s) > 1/2$, as claimed in
Theorem 1.2.

\vskip .10in \flushpar {\bf Meromorphic continuation of
$L_{m}^{+}$: }  The argument to prove the meromorphic
continuation of $L^{+}_{m}$ is similar, in spirit, to that of
$L_{m}^{-}$.
Recall equation (3.1), which shows that
$$
sV_{\ca m}(z,s-1)=2i\frac{d}{dz}U_{\ca m}(z,s)+4\pi mV_{\ca
m}(z,s)\tag9.9
$$
for $\Re (s)$ sufficiently large.  By comparing the series $V_{\ca m}(z,s)$ with $E_{\ca}(z,s+1)$ we see that it converges absolutely and uniformly to a holomorphic function of $s$ for $\Re (s)>0$.
The techniques of Lemma 8.1
apply to  $V_{\ca m}$ to give, for $\Re (s)>0$,
$$
y V_{\ca m}(z,s)\ll 1.
$$
Combining this with (8.18) and (9.9) easily shows that the analytic continuation
of $yV_{\ca m}(z,s-1)$ down to $\Re (s) >1/2$ is bounded by a polynomial in $y_\G (z)$.  Therefore the inner product
$\langle f\K_1, V_m(\cdot,\overline{s}-1)\rangle_2$ admits a
meromorphic continuation for $\Re (s) > 1/2$, which is holomorphic
at $s=1$.  By Proposition 5.1, this implies the meromorphic
continuation of $L^{++}_{m}$ to $\Re (s) > 1/2$.  Since $L^{++}_{m}$ and $L^{+}_{m}$
differ by a Dirichlet polynomial, this part of Theorem 1.2 is now
complete. 

\vskip .25in \head \bf \S 10. Examples
\endhead

\vskip .10in To conclude this work, we will remind the reader of
certain known computations as well as pose a question that can
lead to future investigations.

\vskip .10in Let us consider the discrete subgroup $\PSL_2(\Z)$.
In this case, the Fourier expansion of the first-order Kronecker
limit function is well-known, namely
$$
\K_1(z)=\sum_{n<0}k(n) e(n\overline{z}) + y+K-\frac{3}{\pi}\log y
+ \sum_{n>0}k(n)e(nz) \tag 10.1
$$
where $K=\frac{3}{\pi}(\g - \log 4\pi)$, $\sigma(n)=\sum_{d|n} d$
and
$$
k(n)= \frac{6}{\pi}\frac{\sigma(|n|)}{|n|}.
$$
Also, let us set the notation that for $l \geqslant 0$, we define the
function $\sigma_l(n)=\sum_{d|n} d^l$.  Now consider the
congruence subgroup $\G_{0}(N)$, and, for simplicity, assume that
$N$ is square-free.  As stated in \cite{C-I}, one can express the
first-order non-holomorphic Eisenstein series on $\G_{0}(N)$
through the formula
$$
E(z,s)_{\G_0(N)}=\zeta_N(2s) \sum_{d| N} \mu(d)
(dN)^{-s}E(Nz/d,s),
$$
where $\zeta_N(s)$ is the incomplete zeta-function
$$
\zeta_N(s)=\prod_{p | N}(1-p^{-s})^{-1}
$$
where the product is over all primes $p$ dividing $N$, $\mu$ is the M\"obius function and $E(z,s)$
denotes the $\G=\PSL_2(\Z)$ Eisenstein series.  In effect, this
formula is a consequence of the Artin formalism associated to the
spectral theory on the quotient space $\G_{0}(N)\backslash \H$
viewed as a finite degree cover of $\PSL_{2}(\Z)$.  In the special
case when $N$ is equal to a prime, which we denote by $p$, then we
have that
$$
E(z,s)_{\G_0(p)}=\frac{1}{1-p^{-2s}}\left(p^{-s}E(pz,s)-p^{-2s}E(z,s)\right).
$$
Recall that the volume of $\G_0(p)\backslash \H$ is $p+1$ times
the volume of $\PSL_{2}(\Z)\backslash \H$.  Therefore, one can
compute the first-order Kronecker limit function on $\G_{0}(N)$ to
be
$$
\K_1(z)_{\G_0(p)}=\frac{1}{p^2-1}\left(p \K_1(pz)-\K_1(z)\right)
$$
for $p$ prime.  Therefore, for any prime level $p$, we have, in
effect, computed the Fourier coefficients of the second-order
Kronecker limit function in terms of the divisor sums and the
Fourier coefficients of the chosen degree two form $f \in
S_{2}(\G_{0}(N))$.  Of course, the computations required to
extract the special values $L_{m}^{+}(1)$ and $L_{m}^{-}(1)$,
which require analytic continuation, could be formidable.

\vskip .10in For general Fuchsian groups, the first-order
Kronecker limit function $\K_{1}$ is studied in \cite{Gn}.  The
analogue of the Dedekind $\eta$ function and Dedekind sums are also
studied there.  We refer the interested reader to \cite{Gn} for
additional information.

\vskip .10in Finally, we now highlight a question that arises from
Theorem 1.3.  Given a Fuchsian group $\G$ of the first kind and a
parabolic subgroup, one then has a first-order Kronecker limit
function $\K_{1}$.  With this, consider the map from $S_{k}(\G)$
to itself given by
$$
f \mapsto \Pi_{hol}(f\K_1).\tag 10.2
$$
Are there any interesting characteristics of this map which can
then lead to further simplifications in Theorem 1.3?  Consider the
special case when $\G = \PSL_{2}(\Z)$ and $k = 24$.  In this
setting, we will examine two different holomorphic forms.  The
Dedekind delta function
$$
\Delta (z) = e(z)\prod\limits_{n=1}^{\infty}\left(1 -
e(nz)\right)^{24}
$$
is a weight twelve holomorphic form, as is the Eisenstein series
$$
G_{12}(z)=-B_{12}/24+\sum_{n=1}^\infty \sigma_{11}(n) e(nz),
$$
with $-B_{12}/24=691/65520$.  The vector space
$S_{24}(\PSL_{2}(\Z))$ is two-dimensional with basis $\Delta^2$,
$\Delta G_{12}$, \cite{Za1}.  The analogue of the inner product
formula (3.3) for weight $24$ forms is the identity
$$
\Pi_{hol}(\varphi)=\sum_{m=1}^\infty d_m e(mz) \text{ \ \ for \ \
} d_m=\frac{(4\pi m)^{23}}{22!} \langle \varphi,P_{m}(\cdot)_{24}
\rangle_{24}
$$
(see (8.1)).   With this high weight there will be no problem with
the convergence of the Poincar\'e series.  In general, let
$f(z)=\sum_{n>0} a_n e(nz) \in S_{24}(\G)$ and let
$\Pi_{hol}(f\K_1)=\sum_{m>0} d_m e(mz)$.  Then when using (10.1),
we can compute, as in the beginning of section 6, the formula
$$
\pi d_m=6\sum_{l=1}^m \frac{a_l
\sigma(m-l)}{m-l}+6m^{23}\sum_{l=m+1}^\infty \frac{a_l
\sigma(l-m)}{l^{23}(l-m)}+\frac{23 a_m}{4 m}+3 a_m(2\g +\log m
-H_{22});
$$
$H_n$ denotes the harmonic number $1+\frac 12 +\frac 13+ \cdots
+\frac 1n$, and we have used the formula
$$
\int_0^\infty y^n \cdot\log y\cdot e^{-y} \, dy = n!(H_n -\g)
$$
which holds for $n\geqslant 0$.  This general formula allows for
precise numerical computations.  Specifically, we have computed
that
$$
\Pi_{hol}(\Delta^2 \K_1) \approx -0.852857 \Delta^2 + 0.0000214526
\Delta G_{12} \tag 10.3
$$
and
$$
\Pi_{hol}(\Delta G_{12} \K_1) \approx 0.220305 \Delta^2 +
-0.591762 \Delta G_{12},\tag 10.4
$$
and these computations are correct to the number of decimal places
shown.  In conclusion, these computations suggest that the the
linear map $S_{24}(\PSL_{2}(\Z))\rightarrow S_{24}(\PSL_{2}(\Z))$
given by $f\mapsto \Pi_{hol}(f \K_1)$ to be neither zero nor
diagonal.

\vskip .10in At this time, a host of natural questions arise. For
example, given a Fuchsian group $\G$ and a parabolic subgroup, is
the map (10.2) diagonalizable?  If so, then is there a natural
basis of $S_{k}(\G)$ such that the map (10.2) is diagonal?  Is
there any numerical significance to the coefficients in (10.3) and
(10.4)?  This issues certain warrant future investigations.

\vskip .50in
\noindent
Jay Jorgenson \newline
Department of Mathematics \newline
City College of New York \newline
Convent Avenue at 138th Street \newline
New York, NY 10031 \newline
e-mail: jjorgenson$\@$mindspring.com

\vskip .15in \noindent Cormac O'Sullivan \newline Department of
Mathematics \newline Bronx Community College \newline University
Avenue and West 181st Street \newline Bronx, NY 10453 \newline
e-mail: cormac12$\@$juno.com

\end{document}